\documentstyle{amlts}
\begin{document}
\annalsline{157}{2003}
\received{April 15, 1999}
\revised{October 5, 2000}
\startingpage{715}
\def\bye{\end{document}}
 \font\tenrm=cmr10
\def\ritem#1{\item[{\rm #1}]}
\input boxedeps.tex 
\SetepsfEPSFSpecial 
\HideDisplacementBoxes
\def\figin#1#2{
$$
 {\BoxedEPSF{#1 scaled
#2}%
}%
$$
\noindent}
\catcode`\@=11
\font\twelvemsb=msbm10 scaled 1100
\font\tenmsb=msbm10
\font\ninemsb=msbm10 scaled 800
\newfam\msbfam
\textfont\msbfam=\twelvemsb  \scriptfont\msbfam=\ninemsb
  \scriptscriptfont\msbfam=\ninemsb
\def\msb@{\hexnumber@\msbfam}
\def\Bbb{\relax\ifmmode\let\next\Bbb@\else
 \def\next{\errmessage{Use \string\Bbb\space only in math
mode}}\fi\next}
\def\Bbb@#1{{\Bbb@@{#1}}}
\def\Bbb@@#1{\fam\msbfam#1}
\catcode`\@=12

 \catcode`\@=11
\font\twelveeuf=eufm10 scaled 1100
\font\teneuf=eufm10
\font\nineeuf=eufm7 scaled 1100
\newfam\euffam
\textfont\euffam=\twelveeuf  \scriptfont\euffam=\teneuf
  \scriptscriptfont\euffam=\nineeuf
\def\euf@{\hexnumber@\euffam}
\def\frak{\relax\ifmmode\let\next\frak@\else
 \def\next{\errmessage{Use \string\frak\space only in math
mode}}\fi\next}
\def\frak@#1{{\frak@@{#1}}}
\def\frak@@#1{\fam\euffam#1}
\catcode`\@=12

\title{Convex integration for Lipschitz mappings\\
and counterexamples to regularity} 
\shorttitle{Convex integration for Lipschitz mappings} 

 \acknowledgements{The first named author was supported by  a Max Planck Research Award.   The second
named author was supported by grant DMS-9877055 from the NSF and by a Max Planck Research Award.}
 \twoauthors{S.~M\"uller}{V.~\v{S}ver\'{a}k}
\institutions{Max Planck Institute for Mathematics in the Sciences, Leipzig, Germany\\
{\eightpoint {\it E-mail address\/}: sm@mis.mpg.de}\\
\vglue5pt University of Minnesota, Minneapolis, MN\\
{\eightpoint {\it E-mail address\/}: sverak@math.umn.edu
}}

 \newcommand{\mmn}{{M^{m\times n}}}
\newcommand{\NN}{{\bf N}}
\newcommand{\OO}{{{\cal    O}}}
\newcommand{\RR}{{\bf R}}
\newcommand{\Si}{{\Sigma}}
\newcommand{\MM}{{{\cal    M}}}
\newcommand{\FF}{{{\cal    F}}}
\newcommand{\PP}{{{\cal    P}}}
\newcommand{\CC}{{C}}
\newcommand{\Mrc}{{\PP^{{\rm   rc}}}}
\newcommand{\LL}{{{\cal    L}}}
\newcommand{\Krc}{{K^{{\rm   rc}}}}
\newcommand{\Kjrc}{{K_1^{{\rm   rc}}}}
\newcommand{\Krcs}{{K^{{\rm   rc},\Si}}}
\newcommand{\rof}{{R_\OO f}}
\newcommand{\lo}{{\LL(\OO)}}
\newcommand{\tK}{{\tilde K}}
\newcommand{\fo}{{f_0}}
\newcommand{\slr}{{SL(r,\RR)}}
\newcommand{\ve}{{\varepsilon}}
\newcommand{\phie}{{\phi_{\ve}}}
\newcommand{\psie}{{\psi_{\ve}}}
\newcommand{\dist}{{{\rm   dist\,}}}
\newcommand{\mdd}{{M^{2\times 2}}}
\newcommand{\mcd}{{M^{4\times 2}}}
\newcommand{\sdd}{{S^{2\times 2}}}
\newcommand{\Om}{{\Omega}}
\newcommand{\rd}{{\RR^2}}
\newcommand{\UU}{{{\cal    U}}}
\newcommand{\KK}{{{\cal    K}}}
\newcommand{\varkappa}{{\kappa}}
\newcommand{\ints}{{\int_{\sdd}}}
\newcommand{\om}{{\omega}}
\newcommand{\ome}{{\omega_{\varepsilon}}}
\newcommand{\dper}{{\nabla^\perp}}
\newcommand{\itd}{{\int_{\Td}}}
\newcommand{\vf}{{\varphi}}
\newcommand{\tf}{{\tilde f}}
\newcommand{\sym}{{{\rm   sym}}}
\newcommand{\as}{{{\rm   asym}}}
\newcommand{\Td}{{{\mathbf T}^2}}
\newcommand{\diw}{{{\rm   div\,}}}
\newcommand{\rank}{{{\rm   rank\,}}}
\newcommand{\TT}{{\bf T}}
\newcommand{\snn}{{S^{n\times n}}}
\newcommand{\fca}{a fine $C^0$-approximation }
\newcommand{\orc}{{\OO^{{\rm   rc}}}}
\newcommand{\meas}{{{\rm   meas\,}}}
\newcommand{\aabij}{{a^{\alpha\beta}_{ij}}}

\section{Introduction}\label{intro}

In this paper we study Lipschitz solutions of 
partial differential relations of the form 
\begin{equation}\label{eq: main}
\nabla u(x)\in K \qquad\mbox{a.e.\ in $\Om$,}
\end{equation}
where $u$ is a (Lipschitz) mapping of an open set
$\Omega\subset\RR^n$ into $\RR^m$, $\nabla u(x)$ is its gradient
(i.e.\  the matrix 
${\partial u_i(x)}/{\partial x_j},\,{1\le i\le m,\,1\le j\le n}$,
 defined for almost every $x\in\Om$), and $K$ is a subset of the
set $\mmn$ of all real $m\times n$ matrices.
In addition to relation~(\ref{eq: main}), boundary conditions
and other conditions on $u$ will also be considered.

Relation~(\ref{eq: main}) is a special case of
partial differential relations which have been extensively
studied in connection with certain geometrical problems,
such as isometric immersions. For example, the celebrated results
of Nash \cite{Na 54} and Kuiper \cite{Ku 55} and their  far-reaching generalizations
by Gromov \cite{Gr 86} showed striking and completely unexpected
features of the behavior of\break  $C^1$-isometric
immersions of $\RR^n$ to $\RR^{n+1}$, and Lipschitz isometric immersions
of $\RR^n$ to $\RR^n$. A general result describing a large class
of Lipschitz solutions of partial differential relations more general
than~(\ref{eq: main}) can be found in the  book of Gromov \cite[p.\ 218]{Gr 86}.

More recently,  problems concerning solutions of relations 
of the form~(\ref{eq: main})
have been studied in connection with the characterization of absolute
minimizers of variational integrals describing
the elastic energy of crystals exhibiting
interesting microstructures (\cite{BJ 87}, \cite{CK 88}).
An important observation 
which came from this direction \cite{Ba 90} is that
relation~(\ref{eq: main}) can have highly oscillatory solutions
even when the difference of any two (nonidentical) matrices in $K$ has
rank $\ge 2$. 
This situation, which does occur in some very interesting cases,
is not covered by the theorem of Gromov mentioned above.
In technical terms to be explained below, the reason is that
Gromov's {\it$P$-convex hull} of the set $K$ is again $K$ in that situation.
The main result of this paper, Theorem~\ref{main},  
covers many of these cases and
shows that in the Lipschitz case it seems to be more natural
to work with a different hull, which is defined in terms
of rank-one convex functions, and can be significantly larger
than the $P$-convex hull.

As an application of the theorem we give a solution of a long-standing
problem regarding regularity of weak solutions of elliptic systems.
We construct an example of a variational integral 
$I(u)=\int_{\Om}F(\nabla u)$, where $\Om$ is an open disc in $\RR^2$, $u$ is a mapping
of $\Om$ into $\RR^2$, and $F$ is a smooth, strongly quasi-convex function with 
bounded second derivatives, such that the Euler-Lagrange equation of $I$
has a large class of weak solutions which are Lipschitz but not $C^1$ in any
open subset of $\Om$, and have some other ``wild'' features.
This result should be compared with the well-known result of Evans \cite{Ev 86}
which says that {\it minimizers} of $I$ are smooth  outside 
a closed subset of $\Om$ of measure zero.
Our method also gives new conditions on $F$ which are necessary for regularity.
The conditions are expressed in terms of geometrical properties of the gradient mapping
$X\to DF(X)$.
We expect that the method is applicable to other interesting problems.

Our construction is  quite different from well-known counterexamples
to regularity of solutions of elliptic systems, such as \cite{DG 68}, \cite{GM 68},
or\break \cite{HLN 96}. 
We should emphasize, however, that our method does not apply
when $F$ is convex. 
Very recently we became aware of the work of Scheffer \cite{Sch 74}, 
in which important partial results, including counterexamples,
 related to the regularity problem for the elliptic systems described above
were obtained. It seems that the work was never published in a journal
and has not received the attention it deserves.
 The point of view taken in that paper is implicitly quite
similar to ours and in particular
the $T_4$-configurations discussed in Section~4.2 play an important
role in Scheffer's work. At the same time,  
the new  techniques  we develop
 enable us to answer 
questions which  \cite{Sch 74} left open.

\vglue-8pt
\section{Preliminaries}\label{prelim}
\vglue-3pt

Let us first recall the various notions of convexity related to 
lower-semi\-continuity of variational integrals of the form
$I(u)=\int_{\Om}f(\nabla u)$, where $\Om$ is a bounded domain in $\RR^n$,
$u\colon\Om\to\RR^m$ is a (sufficiently regular) mapping, 
and $f\colon\mmn\to\RR$ is a continuous function defined on the set $\mmn$
of all real $m\times n$ matrices.

A function $f\colon\mmn\to\RR$ is {\it quasi-convex} if
$\int_{\Om} (f(A+\nabla\vf)-f(A))\ge 0$ for each $A\in\mmn$ and each
smooth, compactly supported $\vf\colon\Om\to\RR^m$. This definition
was introduced by Morrey (see e.g.\ \cite{Mo 66}) who also proved that the quasi-convexity
of $f$ is necessary and sufficient for the functional $I$ to be lower-semicontinuous
with respect to the uniform convergence of uniformly Lipschitz functions.
It is also necessary and sufficient  for the weak sequential lower-semicontinuity of $I$
on Sobolev spaces $W^{1,p}(\Om,\RR^m)$, if natural growth conditions
are satisfied; see \cite{Ma 85} and \cite{AF 87}. 
The definition of quasi-convexity is independent of $\Om$, as can be seen
by a simple scaling and covering argument (\cite{Mo 66}). 
In fact, we have the following simple observation made by many authors:

\proclaim{Lemma}\label{per}
Let $\TT^n$ be a flat $n$\/{\rm -}\/dimensional torus. A function $f\colon\mmn\break\to\RR$
is quasi\/{\rm -}\/convex if and only if $\int_{\TT^n}(f(A+\nabla\vf)-f(A))\ge0$ for each
$A\in\mmn$ and each smooth $\vf\colon\TT^n\to\RR^m$.
\endproclaim

The reader is referred  to \cite{Sv 92a} for a proof of this statement.

We also recall that, with the notation above, 
$f\colon\mmn\to\RR$ is {\it strongly
quasi-convex} if there exists $\gamma>0$ such that
$\int_{\Om} (f(A+\nabla\vf)-f(A))\ge \gamma\int_{\Om}|\nabla\vf|^2$ 
for each $A\in\mmn$ and each
smooth, compactly supported $\vf\colon\Om\to\RR^m$.
This notion appears naturally in the regularity theory; see for example \cite{Ev 86}.

A function $f\colon\mmn\to\RR$ is {\it rank-one convex} if it is convex
along any line whose direction is given by a matrix of rank one, i.e.\ 
$t\to f(A+tB)$ is convex for each $A\in\mmn$ and each $B\in\mmn$ with
$\rank B=1$. This class of functions will play a particularly important
r\^ole in our analysis. It can be proved that any quasi-convex function is 
rank-one convex, but the opposite implication  fails when $n\ge2,\,m\ge3$ (\cite{Sv 92a}).
(The case $n\ge2,\,m=2$ is open.)

We will also deal with functions which are defined only on symmetric matrices.
We will denote by $\snn$ the set of all symmetric $n\times n$ matrices.
The notions introduced above for functions on $\mmn$ can be modified in the
obvious manner to apply to functions on symmetric matrices.
For example,
a function $f\colon\snn\to\RR$ is quasi-convex, if
$\int_{\Om} (f(A+\nabla^2\phi)-f(A))\ge 0$ 
for each $A\in\snn$ and each
smooth, compactly supported $\phi\colon\Om\to\RR$.
Again, the definition is independent of $\Om$ and, in fact, $\Om$
can be replaced by any flat $n$-dimensional torus.

In the rest of this section we examine in more detail   
facts related to rank-one convexity.

Let $\OO\subset\mmn$ be an open set and let $f\colon\OO\to\RR$
be a function. We say that $f$ is {\it rank-one convex in $\OO$},
if $f$ is convex on each rank-one segment contained in $\OO$.
It is easy to see that every rank-one convex function
$f\colon\OO\to\RR$ is locally Lipschitz in $\OO$.

We will use $\PP$ to denote the set of all compactly supported 
probability measures
in $\mmn$. For a compact set $K\subset\mmn$ we use $\PP(K)$ to denote
the set of all probability measures supported in $K$.
For $\nu\in\PP$ we denote by $\bar\nu$ the center of mass of $\nu$,
i.e.\  $\bar\nu=\int_{\mmn}Xd\nu(X)$.

Following \cite{Pe 93}, we say that a measure $\nu\in\PP$ is a {\it laminate} if 
$\langle\nu,f\rangle\ge f(\bar\nu)$ for each rank-one convex
function $f\colon\mmn\to\RR$. 
At the center of our attention will be the
sets $\Mrc(K)=\{\nu\in\PP(K),\  \nu\ {\rm   is\ a\ laminate}\}$,
which are defined for any compact set $K\subset\mmn$.

For $A\in\mmn$ we denote by $\delta_A$ the Dirac mass at $A$.

Let $\OO$ be an open subset of $\mmn$. Assume $\nu\in\PP$ is of the
form $\nu=\sum_{j=1}^{j=r}\lambda_j\delta_{A_j}$, with 
$A_j\in\OO$, $j=1,\dots,r$, and $A_j\ne A_k$ when $j\ne k$. 
We say that $\nu'\in\PP$ can be obtained
from $\nu$ by an {\it elementary splitting in $\OO$} if, for some
$j\in\{1,\dots,r\}$, and some $\lambda\in[0,1]$,
 there exists a rank-one segment 
$[B_1,B_2]\subset\OO$ containing $A_j$, with $A_j=(1-s)B_1+sB_2$,
such that 
 $\nu'=\nu+\lambda\lambda_j((1-s)\delta_{B_1}+s\delta_{B_2}-\delta_{A_{j}})$.

We now define an important subset $\LL(\OO)$ of laminates, called
{\it  laminates of a finite order in $\OO$}. By definition, $\nu\in\LL(\OO)$
if there exists a finite sequence of measures 
$\nu_1,\dots,\nu_m$ such that $\nu_1=\delta_A$ for some $A\in \OO$,
$\nu_m=\nu$, and $\nu_{j+1}$ can be obtained from $\nu_j$ by an
elementary splitting in $\OO$ for $j=1,\dots,m-1$. When $\OO=\mmn$, the measures in
$\LL(\OO)=\LL(\mmn)$ are called {\it laminates of a finite order} (i.e.\  we 
do not refer to the set $\OO$ in that case).

 Let $K$ be a compact subset of $\mmn$. 
The {\it rank-one convex hull $\Krc\subset\mmn$ of $K$}
  is defined as follows. A matrix $X$ does not belong to
$\Krc$ if and only if there exists
$f\colon\mmn\to\RR$ which is 
rank-one convex  such that $f\le 0$ on $K$
and $f(X)>0$. We emphasize that this definition will be used only
when $K$ is compact. For open sets $\OO\subset\mmn$, we define
the rank-one convex hull $\orc$  of $\OO$ as $\orc=\mathbold{\cup}\{\Krc,\,\,
\mbox{$K$ is a compact subset of $\OO$}\}$.
With this definition we have the property that the rank-one convex
hull of an open set is again an open set, which will be useful for
our purposes. 

We refer the reader to \cite{MP 98} for interesting results about
rank-one convex hulls of closed sets. 
The following theorem, which is a slight generalization of a result from \cite{Pe 93}, will 
play an important r\^ole.  
\def\spn#1{\specialnumber{#1}}
\specialnumber{2.1}\proclaim{Theorem}\label{lo-thm}
Let $K$ be a compact subset of $\mmn$ and let 
$\nu\in\Mrc(K)$. 
Let $\OO\subset\mmn$ be an open set such that $\Krc\subset\OO$.
Then there exists a sequence $\nu_j\in\LL(\OO)$ of laminates of a finite order
in $\OO$ such that $\bar\nu_j=\bar\nu$ for each $j$ and the
$\nu_j$ converge weakly$^*$ to $\nu$ in $\PP$.
\endproclaim

As a preparation for the proof of the theorem, we prove the following 
lemma.
\spn{2.2}\proclaim{Lemma}\label{l1}
Let $\OO$ be an open subset of $\mmn$. 
Let $f\colon\OO\to\RR$ be a continuous function and let
$\rof\colon\OO\to\RR\cup\{-\infty\}$ be defined by
$$
\rof=\sup\{\mbox{$g,\, g\colon\OO\to\RR$ is rank\/{\rm -}\/one convex in $\OO$
and $g \le f$} \}.
$$
Then for each $X\in\OO${\rm ,}
 $\rof(X)=\inf\{\mbox{$\langle\nu,f\rangle,\, \nu\in\LL(\OO)$ and
$\bar\nu=X$}\}$.
\endproclaim
\demo{{P}roof}
Let us denote by $\tilde f$ the function in $\OO$ defined by
$\tilde f(X)=\inf\{\langle\nu,f\rangle,$\break $\nu\in\LL(\OO)\ {\rm   and\ }
\bar\nu=X\}$. Clearly $\rof\le\tilde f$ in $\OO$. On the other hand,
we see from the definition of the set $\LL(\OO)$ that it has the
following property: if $\nu_1, \nu_2\in\lo$,  and the segment
$[\bar\nu_1,\bar\nu_2]$ is a rank-one segment contained in $\OO$, then any convex
combination of $\nu_1$ and $\nu_2$ is again in $\LL(\OO)$. Using this,
we see immediately from the definitions that $\tilde f$ is
rank-one convex in $\OO$ and hence $\rof=\tilde f$.
\enddemo

{\it Proof of Theorem}~\ref{lo-thm}. 
Let $\nu\in\Mrc(K)$ and let $\bar\nu=A$ be its center of mass.
From the definitions we see that $A\in\Krc$. 
We choose an open set $U\subset\mmn$  satisfying
$\Krc\subset U\subset\bar U\subset \OO$ and define
$\FF=\{\mu\in\LL(U),\ \bar\mu=A\}$. We claim that the weak$^*$ closure
of $\FF$ contains $\nu$. To prove the claim, we argue by contradiction.
Assume $\nu$ does not belong to the weak$^*$ closure of $\FF$.
Since $\FF$ is clearly convex, we see from the Hahn-Banach theorem
that there exists a continuous function $f\colon\bar U\to\RR$ such that
$\langle\nu,f\rangle<\inf\{\langle\mu,f\rangle,\ \mu\in\LL(U)$ and $
\bar\mu=A\}$. By Lemma~\ref{l1}, we have
$\inf\{\langle\mu,f\rangle,\ \mu\in\LL(U)\ {\rm   and}\
\bar\mu=A\}=R_Uf(A)$. We see that the function $\tilde f=R_Uf\colon
U\to\RR$ is rank-one convex in $U$ and satisfies 
$\langle\nu,\tilde f\rangle\le
\langle\nu, f\rangle < \tilde f(\bar\nu)$. By Lemma~\ref{ext-lemma} below,
there exists a rank-one convex function
$F\colon\mmn\to\RR$ such that $F=\tilde f$ on $\Krc$. We conclude
that $\nu$ cannot belong to $\Mrc(K)$, a contradiction.
The proof is finished.

\spn{2.3}\proclaim{Lemma}\label{ext-lemma}
Let $K\subset\mmn$ be a compact set, let $\OO$ be an open set containing
$\Krc$ {\rm (}\/the rank\/{\rm -}\/one convex hull of $K$\/{\rm )}
 and let $f\colon\OO\to\RR$ be rank\/{\rm -}\/one convex. Then there exists
$F\colon\mmn\to\RR$ which is rank\/{\rm -}\/one convex and coincides with $f$ in a
neighborhood of $\Krc$.
\endproclaim

\demo{{P}roof} We claim there exists a nonnegative rank-one convex 
$g\colon\mmn\break\to\RR$ such that $\Krc=\{X,\,g(X)=0\}$.
 To prove this, we choose $R>0$ so that $K\subset B_{R/2}=\{X, |X|<R/2\}$
and define $g_1\colon B_R\to\RR$ by 
\begin{eqnarray*}
 g_1(X)&=&\sup\{f(X),\,f\colon B_R\to\RR, \\
 & &\quad \mbox{$f$ is
 rank-one convex in $B_R$ and $f\le\dist(\,\cdot\,,K)$ in $B_R$}\}.
\end{eqnarray*}
The function $g_1$  is obviously nonnegative and rank-one convex
in $B_R$. Moreover, 
 $\{X\in B_R,\,g_1(X)=0\}\supset K$ and from the definition of $\Krc$ we see that
$g_1>0$ outside $\Krc$. We now define
$$ g(X)=\left\{ \begin{array}{ll}
                  \max\,(g_1(X), 12|X|-9R) & \mbox{when $X\in B_R$} \\
               12|X|-9R &  \mbox{when $|X|\ge R$.}
                 \end{array} \right. $$
Clearly $g$ is rank-one convex  in a neighborhood of any point
$X$ with $|X|\ne R$. Since $g_1(X)\le 2|X|$ when $|X|=R$, we see
that we have $g(X)=12|X|- 9 R$ in a neighborhood of $\{|X|=R\}$.
We see that $g$ is nonnegative, rank-one convex in $\mmn$, $\{X, g(X)=0\}\supset K$, and
$\{X, g(X)>0\}\cap\Krc=\emptyset$. Therefore $\{X, g(X)=0\}=\Krc$

We can now finish the proof of the lemma. 
 Replacing $f$ by $f+c$, if necessary,
we can assume that $f>0$ in a neighborhood of $\Krc$.
For $k>0$ we let $U_k=\{X\in\OO,\,\,f(X)>kg(X)\}$.
We also let $V_k$ be the union of the connected components
of $U_k$ which have a nonempty intersection with $\Krc$.
It is easy to see that there exists $k_0>0$ such that
$\bar V_{k_0}\subset\OO$. We now let $F(X)=f(X)$ when $X\in V_{k_0}$
and $F(X)=k_0g(X)$ when $X\in\mmn\setminus V_{k_0}$. 
It is easy to check that the function $F$
defined in this way is rank-one convex\break on $\mmn$. 

\section{Constructions}\label{constr}

Throughout this section, $\Om$ denotes a fixed bounded open subset of $\RR^n$.
We will use the following terminology. A Lipschitz 
mapping $u\colon\Om\to\RR^m$ is {\it piecewise affine},  if there exists
a countable system of mutually disjoint open sets $\Om_j\subset\Om$ 
which cover $\Om$ up to a set of zero measure, and the restriction of $u$ to
each of the sets $\Om_j$ is affine.

Following Gromov (\cite[p.\ 18]{Gr 86}) we also introduce the following concept.
Let $\FF(\Om,\RR^m)$ be a family of continuous mappings of $\Om$ into $\RR^m$.
We say that a given continuous mapping $v_0\colon\Om\to\RR^m$ admits
{\it a fine $C^0$-approximation} by the family $\FF(\Om,\RR^m)$ if there exists,
for every continuous 
function $\ve\colon\Om\to(0,\infty)$, an element $v$ of the family
$\FF(\Om,\RR^m)$ such that $|v(x)-v_0(x)|<\ve(x)$ for each $x\in\Om$.

\demo{{\rm 3.1.} The basic construction}
The main building block of all the solutions of relation~(\ref{eq: main}) which 
we construct in this paper is the following simple lemma.
\enddemo

\proclaim{Lemma}\label{basic}
Let $A,B\in\mmn$ be two matrices with $\rank (B-A)=1${\rm ,} let $b\in\RR^m${\rm ,}
$0<\lambda<1$ and  $C=(1-\lambda)A + \lambda B$. 
Then{\rm ,} for any $0<\delta<|A-B|/2${\rm ,}
 the affine mapping $x\to Cx+b$ admits  \fca 
by piecewise affine mappings
$u\colon\Om\to\RR^m$
such that $\dist(\nabla u(x),\{A,B\})<\delta$
almost everywhere in $\Om${\rm ,} $\meas\{x\in\Om,\,\,
|\nabla u(x)-A|<\delta\} = (1-\lambda)\,\meas\Om${\rm ,}
and 
$\meas\{x\in\Om,\,\,
|\nabla u(x)-B|<\delta\} = \lambda\,\meas\Om$.
\endproclaim 
 
{\it Proof}. We first note that it is enough to prove the lemma only for a
special case when the function $\ve(x)$ appearing in the definition of a fine\break $C^0$-approximation is 
constant and the function approximating the function $u$ satisfies the boundary
condition $u(x)=Cx+b$ for $x\in\partial\Om$.
 This can be seen by considering a sequence of open sets $\Om_j$ which
are mutually disjoint, satisfy $\bar\Om_j\subset\Om$, and cover $\Om$ up
to a set of measure zero.

To prove the special case, we note that we can assume without loss of generality
that $A=-\lambda a\otimes e_n,\ B=(1-\lambda)a\otimes e_n$, and $C=0$,
where $a\in\RR^m$ and $e_n=(0,\dots,0,1)\in\RR^n$.
We define $h\colon\RR\to\RR$ and $w\colon\RR^n\to\RR^m$ by
$h(s)=(|s|+(2\lambda-1)s)/2$ and $w(x)=a\max(0,1-|x_1|-\dots-|x_{n-1}|-h(x_n))$.
We choose a small $\delta'>0$, and set
$v(x)=\delta' w(x_1,\dots,x_{n-1},x_n/\delta')$. We also let $\omega=\{x,\, v(x)>0\}$.
We check by a direct calculation that $\dist(\nabla v(x),\{A,B\})\le (n-1)|a|\delta'$
for almost every $x\in\omega$. We clearly also have $v(x)=0$ when $x\in\partial\omega$.
By Vitali's theorem we can cover $\Omega$ up to a set of measure zero
 by a countable family $\{\omega_i\}$ of mutually
disjoint sets of the form $\omega_i=y_i+r_i\omega$ (with $y_i\in\RR^n$ and 
$r_i\in(0,\epsilon)$). We let
$u(x)=r_iv(r_i^{-1}(x-y_i)$ when $x\in\omega_i$, and $u(x)=0$ if 
$x\in\Om\setminus\cup_i\omega_i$. It easy to check that $u$ satisfies the
required conditions, provided $\delta'$ is sufficiently small.

\proclaim{Lemma}\label{sl}
Let $\nu\in\PP(\mmn)$ be a laminate of a finite order{\rm ,} let $A=\bar\nu$
be its center of mass. 
Let us write $\nu=\sum_{j=1}^r\lambda_j\delta_{A_j}$ with
$\lambda_j>0$ and $A_i\ne A_j$ when $i\ne j${\rm ,} and let $$\delta_1=\min\{
|A_i-A_j|/2; 1\le i<j\le r\}.$$
Then{\rm ,} for each $b\in\RR^m$, and each $0 < \delta < \delta_1${\rm ,}
the mapping $x\to Ax \! + b$\break admits \fca by piecewise affine mappings $u$
satisfying\break $\dist(\nabla u(x),\{A_1,\dots,A_r\})<\delta$ {\rm a.e.}\
in $\Om$ and $$\meas\{x\in\Om,\, \dist(\nabla u(x), A_j) <\delta\}=\lambda_j\,
\meas\Om$$ for each $j\in\{1,\dots,r\}$.
\endproclaim

{\it Proof}. This can be easily proved by applying iteratively Lemma~\ref{basic}
 in a way which is naturally suggested by
 the definition of the laminate
of a finite order. We outline some details for the convenience of the reader.
Let $\delta_A=\nu_1,\nu_2,\dots,\nu_m=\nu$ be a sequence of measures
such that $\nu_{j+1}$ can be obtained from $\nu_j$ by an elementary
splitting in $\mmn$. If $m=1$, there is nothing to prove, if $m=2$,
our statement is exactly Lemma~\ref{basic}. Proceeding by induction
on $m$, let us assume that the lemma has been proved for $\nu$ replaced
by $\nu_{m-1}$. Let us write $\nu_{m-1}=\sum_{j=1}^{j=r'}\lambda'_j
\delta_{A'_j}$, with $A'_k\ne A'_l$ when $k\ne l$. Since $\nu=\nu_m$
can be obtained from $\nu_{m-1}$ by an elementary splitting, 
$$\nu=\nu_{m-1}+\lambda\lambda'_{j_0}
((1-s)\delta_{B_1}+s\delta_{B_2}-\delta_{A'_{j_0}})$$
for some $\lambda\in[0,1]$, $s\in[0,1]$, $j_0\in\{1,\dots,r'\}$, and
a rank-one segment $[B_1,B_2]$ containing $A'_{j_0}$. By our assumptions, for any sufficiently
small $0<\delta'<\delta/2$, the map $x\to Ax+b$ admits \fca
by piecewise affine maps $u'$ satisfying
$\dist(\nabla u'(x),\{A'_1,\dots,A'_{r'}\})<\delta'$ a.e.\  in $\Om$
and $$\meas\{x\in\Om; \dist(\nabla u'(x), A'_j)<\delta'\}=\lambda_j'\,\meas\Om.$$
For any such $u'$ we can find an open set $\Om'\subset\Omega$ such that 
$\dist(\nabla u'(x),A'_{j_0})<\delta'$ in $\Om'$,
$\meas\Om'=\lambda\,\meas\{x\in\Om; \dist(\nabla u'(x),A'_{j_0})<\delta'\}=
\lambda\lambda'_{j_0}\,\meas\Om$, and  $u'$ is piecewise affine in 
$\Omega'$.  
Let $\Om'_k\subset \Omega', \,\,k=1,2,\dots$ 
be mutually disjoint open sets which cover $\Om'$ up to
a set of \pagebreak measure zero such that $\nabla u'=\tilde A_k={\rm const}$ 
in $\Om'_k$, with $|\tilde A_k-A'_{j_0}|<\delta'$.
We now adjust $u'$ by applying Lemma~\ref{basic} on each $\Om'_k$ with
$A=B_1+\tilde A_k-A'_{j_0}$, $B= B_2+\tilde A_k-A'_{j_0}$, 
$C=\tilde A_k$, $\delta=\delta'$, and the proof is easily finished.

\vglue8pt

3.2. {\it Open relations}.
We recall that the rank-one
convex hull $\orc$ of an open set $\OO\subset\mmn$ is, by definition, 
the union of the rank-one
convex hulls of all compact subsets of $\OO$. The main result of this
subsection is the following.
 
\spn{3.1}\proclaim{Theorem}\label{open}
Let $\OO\subset\mmn$ be open{\rm ,} and let $P\subset\orc$ be compact.
Let $u_0\colon\Om\to\RR^m$
be a piecewise affine Lipschitz mapping such that
$\nabla u_0(x)\in P$ for {\rm a.e.}\  $x\in\Om$.
Then $u_0$ admits \fca by piecewise affine
Lipschitz mappings $u\colon\Om\to\RR^m$
satisfying $\nabla u(x)\in\OO$ {\rm a.e.}\  in $\Om$.
\endproclaim
 
{\it Proof}.  As a first step, we prove the following lemma.
 
\spn{3.3}\proclaim{Lemma}\label{stepone} \hskip-8pt
Let $K\subset\mmn$ be a compact set and let $U\subset\mmn$ be an open
set containing $K$. Let $\nu\in\Mrc(K)$ and denote $A=\bar\nu$. Let $b\in~\RR^m$.
Then{\rm ,} for any given $\delta>0${\rm ,} the mapping $x\to Ax+b$ admits
a fine\break $C^0$\/{\rm -}\/approximation by piecewise affine mappings
$u$ satisfying $\nabla u(x)\in U^{{\rm   rc}}$ {\rm a.e.}\
in $\Omega$ and $\meas\{x\in\Om,\,\nabla u(x)\in U\}>(1-\delta)\,\meas\Omega$.
\endproclaim

{\it Proof}. By Theorem~\ref{lo-thm} there exists a laminate $\mu$ of a 
finite order which is supported in a finite subset of 
$U^{{\rm   rc}}$ and satisfies $\bar\mu=\bar\nu$ and $\mu(U)>(1-\delta)$.
Let us write $\mu=\sum_{j=1}^{j=r}\lambda_j\delta_{A_j}$, so that 
$\delta_1=\min\{|A_k-A_l|/2;\,1\le k<l\le r\}>0$. We choose $0<\delta'<\delta_1$
so that each $A_k\in U$ is at distance at least $\delta'$ from the boundary
$\partial U$. From Lemma~\ref{sl} we see that the map $x\to Ax+b$ admits
\fca by piecewise maps $u$ such that
$\dist(\nabla u(x), \{A_1,\dots,A_r\})\break<\delta'$ a.e.\  in $\Omega$
and $\meas\{x\in\Om;\,\dist(\nabla u(x), A_j)<\delta'\}=\lambda_j\,\meas\Om$
for $j=1,\dots,r$, and our lemma immediately follows.
\vglue4pt

Theorem~\ref{open} can now be proved by repeatedly applying  Lemma~\ref{stepone}
in the following way. We first choose a sequence of compact sets
$K_1,K_2,\dots\subset\mmn$, a sequence of open sets $U_1,U_2,\dots\subset\mmn$, 
and a compact set $Q\subset\mmn$ such that
$P=K_1\subset U_1\subset K_2\subset U_2\subset\dots\subset Q\subset \orc$.
We also choose $0<\delta<1$. 
Let $\ve=\ve(x)>0$ be a continuous function on $\Omega$.
In the first step we apply Lemma~\ref{stepone} to approximate $u_0$ up to  $\ve/2$
by a mapping $u_1$ satisfying $\nabla u_1(x)\in U_1^{{\rm   rc}}$ a.e.\  
in~$\Omega$,
together with $\meas\{x\in\Omega,\,\nabla u_1(x)\break\in U_1\}>(1-\delta)\meas\Omega$.
We now modify $u_1$ on those subregions of $\Omega$ where $\nabla u_1(x)$ does not
belong to $U_1$ by applying Lemma~\ref{stepone} again. We obtain a new mapping,
$u_2$, which approximates $u_1$ up to $\ve/4$, coincides with $u_1$ a.e.\ 
in the set $\{x\in\Omega,\,\nabla u_1(x)\in U_1\}$, and satisfies 
$\nabla u_2(x)\in U_2^{{\rm   rc}}$
a.e.\ in $\Omega$ together with
 $\meas\{x\in\Omega,\,\nabla u_2(x)\in U_2\}>((1-\delta)+\delta(1-\delta))\,\meas\Omega$.
By continuing this procedure we get a sequence $u_k$ of mappings which is easily seen to
converge to a mapping
$u$ which gives the required approximation of $u_0$.
\pagebreak

{\it Remark}. From the proofs of Lemma~\ref{sl}, Lemma~\ref{stepone}, and
Theorem~\ref{open} it is easy to see that Lemma~\ref{sl} remains true
if $\nu$ is a laminate (not necessarily of finite order) which can be
written as a finite convex combination of Dirac masses.
\vglue6pt

3.3. {\it Closed relations and in\/{\rm -}\/approximations}.
When considering relation~(\ref{eq: main}) for closed sets $K$, it is natural
to try to construct solutions by combining Theorem~\ref{open} and a suitable
limit procedure. For simplicity we will assume in this section that $K$ is compact.
Following Gromov (\cite[p.\ 218]{Gr 86}) we say that a sequence of open sets
$\{U_i\}_{i=1}^{\infty}$ is an {\it in-approximation} of $K$ if 
$U_i\subset U_{i+1}^{{\rm   rc}}$ for each $i$,
and $\,\,\sup_{X\in U_i}\dist(X,K)\to 0$ as $i\to\infty$. (The definition does not 
require
that each point of $K$  can be reached by a sequence $X_j\in U_j$.)
\spn{3.2}\proclaim{Theorem}\label{main}
Assume that a compact set $K\subset\mmn$ admits an in\/{\rm -}\/approxi\-mation by open
sets $U_i$ in the sense of the definition above. Then any 
 $C^1$\/{\rm -}\/mapping
$v\colon\Omega\to\RR^m$ satisfying $\nabla v(x)\in U_1$ in $\Omega$
admits a fine $C^0$\/{\rm -}\/approxi\-mation 
by Lipschitz mappings $u\colon\Omega\to\RR^m$
satisfying $\nabla u(x)\in K$ \/{\rm a.e.}\/\ in $\Omega$.
\endproclaim

\vglue-12pt
{\it Proof}.
By the same argument as in the proof of Lemma~\ref{basic} it is enough
to prove the statement only in the case when the function $\ve=\ve(x)$ in the
definition of a fine $C^0$-approximation is constant.

Let $\rho\colon\RR^n\to\RR$ be the usual mollifying kernel, i.e.\ we assume that 
$\rho$ is smooth, nonnegative, supported in $\{x,\ |x|<1\}$, and $\int \rho=1$.
For $\ve>0$ we let $\rho_{\ve}=\ve^{-n}\rho(x/\ve)$.
For a function $w\in L^1(\Omega)$ we define $\rho_{\ve} * w$ in the usual way,
by considering $w$ as a function on $\RR^n$ with $w=0$ outside $\Omega$.
In other words, $\rho_{\ve} * w(x)=\int_{\Omega}w(y)\rho_{\ve}(x-y)\,dy$.

We start the proof by choosing $\delta_1>0$ (the exact value of which will be
specified later) and by approximating $v$ by a piecewise affine 
$u_1\colon\Omega\to\RR^m$ with $|u_1-v|<\delta_1$ in $\Omega$, $u_1=v$ on 
$\partial\Omega$, and $\nabla u_1\in U_1$ a.e.\ in $\Omega$.
(We recall that in this paper ``piecewise affine'' allows for countably
many affine pieces.)
We also choose $\ve_1>0$ so that 
$||\nabla u_1 * \rho_{\ve_1}-\nabla u_1||_{L^1(\Omega)}\le 2^{-1}$.

Using Theorem~\ref{open} together with an obvious inductive argument,
we construct a sequence of mappings $u_i\colon\Omega\to\RR^m$ and numbers
$0<\ve_i<2^{-i},\, \delta_i>0$ satisfying
$$\begin{array}{rll}
\noalign{\vskip-4pt}
\nabla u_i & \in\ \,   U_i &\qquad \mbox{a.e.\  in $\Omega$}\,,\\[3pt]
 u_i & =\ \, v &\qquad  \mbox{on $\partial\Omega$}\,,\\[3pt]
||\nabla u_i * \rho_{\ve_i}-\nabla u_i||_{L^1(\Omega)} & \le\ \,  2^{-i}\,,\\[3pt]
\delta_{i+1} & =\  \, \ve_i\delta_i\, ,\\[3pt]
|u_{i+1}-u_i| & \le\ \, \delta_{i+1} &\qquad \mbox{in $\Omega$}\,.
\end{array}
$$
The mappings $u_i$ converge uniformly to a Lipschitz function
$u\colon\Omega\to\RR^m$. We also have
$|u-v|\le\sum_i |u_{i+1}-u_i| + |u_1-v|\le 2\delta_1$.
It remains to  \pagebreak prove that $\nabla u\in K$ a.e. in $\Omega$.
This will be clear if we establish that $\nabla u_i\to\nabla u$
in $L^1(\Omega)$.
We can write 
\begin{eqnarray*}
||\nabla u_i - \nabla u||_{L^1(\Omega)}
 & \le & ||\nabla u_i - \nabla u_i*\rho_{\ve_i}||_{L^1(\Omega)}\\[4pt]
 && +\  ||\nabla u*\rho_{\ve_i} - \nabla u||_{L^1(\Omega)}\\[4pt]
 && +\   ||\nabla u_i*\rho_{\ve_i} - \nabla u*\rho_{\ve_i}||_{L^1(\Omega)}.
\end{eqnarray*}
 
The first two terms on the right-hand side of this inequality clearly converge
to zero as $i\to\infty$. 
 Defining
         $\Omega_i = \{ x \in \Omega, \mbox{dist}(x,\partial \Omega) > 2
         \ve_i \}$ we can  estimate the third term as
      $$
           ||(u_i-u)*\nabla\rho_{\ve_i}||_{L^1(\Omega)}
           + ||\nabla u_i - \nabla u||_{L^1(\Omega \setminus \Omega_i)}
                    \le
          \frac c{\ve_i}||u_i-u||_{\infty} + C \,\meas 
           (\Omega \setminus \Omega_i)\,,
      $$
          where $c$ and $C$ are constants depending only on $\rho$ and
        the Lipschitz constant of $u_i-u$, respectively.

 We have
\begin{eqnarray*}
||u_i-u||_{\infty} &\le& \sum_{j=i}^{\infty}||u_j-u_{j+1}||_{\infty}
 \le \sum_{j=i+1}^{\infty}\delta_j\le 2\delta_{i+1}.
\end{eqnarray*}
Hence the third term can be estimated by 
$$2c\delta_{i+1}/\ve_i+ 
C \,\meas (\Omega \setminus \Omega_i)\le2c\delta_i+ C\, \meas (\Omega \setminus \Omega_i)$$
which converges to zero as $i\to\infty$. The proof is finished.
\vglue8pt

{\it Remark}. The explanation of the strong convergence of $\nabla u_i$
is more or less the following. We can achieve a very fast convergence of
$u_i$ in the sup-norm. It may seem that this is not enough to say much 
about the convergence of $\nabla u_i$. However, in the proof we choose
the parameters in such a way that \mbox{$||u_i-u||_{\infty}$} is very small
in comparison with a typical length over which $\nabla u_i$ changes
significantly (in an integral sense). Therefore, as regards the convergence of
$\nabla u_i$, we get a situation which is in a certain sense
similar to the simple case when the functions $u_i$ are affine in $\Om$. 
This is the main reason we get
the strong convergence. The above argument is taken from \cite{MS 96}.
A different approach can be found in \cite{DM 97}.

\section{Applications to elliptic systems}\label{applications}
\advance\eqcount by 1

Let $\Om\subset\rd$ be a disc. For (sufficiently regular) 
mappings $u\colon\Om\to\rd$ we consider the functional
$I(u)=\int_{\Om}F(\nabla u(x))\,dx$, where $F$ is a (smooth) function
on the set $\mdd$ of all real $2\times 2$ matrices, which satisfies
certain ``ellipticity conditions''. More precisely, we will require
that $F$ be strongly quasiconvex and that its second derivatives be
uniformly bounded in $\mdd$.

The purpose of this section is to show how we can apply the results above
to construct  weak solutions of the Euler-Lagrange equation
\begin{equation}
\diw DF(\nabla u)=0 \label{eq: e-l}
\end{equation}
of the functional $I$ which are Lipschitz, but not continuously
differentiable on any open subset of $\Om$.
This is in sharp contrast with regularity properties of minimizers of $I$,
see, for example \cite{Ev 86}. In fact, we prove the following slightly stronger
statement.

\proclaim{Theorem}\label{solutions}
There exists a smooth strongly quasiconvex function 
$F_0\colon\mdd\to\RR$ with $|D^2F_0|\le c$ in $\mdd${\rm ,}
four matrices $A_1,\dots,A_4\in\mdd${\rm ,} $\varepsilon>0$ and $\delta>0$
such that the following is true.
Let $F\colon\mdd\to\RR$ be a\break $\CC^2$\/{\rm -}\/function satisfying
\hbox{$|DF(A_j)\!-\!DF_0(A_j)|\le\delta\!$ and $|D^2F(A_j)\!-\!D^2F_0(A_j)|\le\delta$}
 for $j=1,2,3,4.$
Then each piecewise 
$\CC^1$\/{\rm -}\/function $v\colon\Omega\to\rd$ satisfying\break $|\nabla v|<\ve$  
{\rm a.e.}\ 
in $\Omega$ admits a fine $\CC^0$\/{\rm -}\/approximation by Lipschitz mappings
$u\colon\Om\to\rd$ which are not $\CC^1$ on any open subset of $\Omega$
and are weak solutions of the equation 
$\diw DF(\nabla u)=0 $ in $\Omega$.
\endproclaim

The theorem will be proved in Section~4.4, after we establish some useful
facts about quasiconvex functions and rank-one convex hulls.
The idea of the construction is the following. We rewrite 
equation~(\ref{eq: e-l}) as a first-order system 
\begin{equation}
\nabla w\in K \label{eq: sys}
\end{equation}
and then show that the strong quasiconvexity does not prevent the rank-one
convex hull of $K$ from being large. (We note that the strong quasi-convexity
does exclude any nontrivial rank-one connections in $K$; see \cite{Ba 80}.)
 We can then use the 
methods developed in the previous sections to construct the desired solutions.
Moreover, it turns out that the situation is stable under the perturbations of $F_0$ which
are allowed in the theorem.
 
\vglue12pt
{\it Remark}. In \cite{Sch 74} Scheffer constructs counterexamples  
to partial regularity of solutions of equation (\ref{eq: e-l}) with
$F$ rank-one convex and with  $u$ in the Sobolev space $W^{1,1}$.
\vglue12pt

One way to write equation~(\ref{eq: e-l}) in the form~(\ref{eq: sys})
is the following. We denote by $J$ the matrix 
$\left(\begin{array}{rr}
         0 & -1\\
         1 & 0
         \end{array}\right)^{\phantom{|}}\!\!$.
The condition that the $2\times 2$ tensor $DF(\nabla u)$ be 
divergence-free is
equivalent to the condition that $DF(\nabla u)J$ be the gradient
of a function $\tilde u\colon\Om\to\RR^2$. We now  \pagebreak introduce
 $w\colon\Om\to\RR^4$ by 
$w=\left(\begin{array}{l} u \\ \tilde u \end{array}\right)^{\phantom{|}}\!\!$. We also let
$K$ be the set of all $4\times 2$ matrices of the form
$\left(\begin{array}{c} X \\ DF(X)J \end{array}\right)$, where
$X$ runs through all $2\times 2$ matrices. It is clear that,
in this notation, system~(\ref{eq: e-l}) is equivalent
to system~(\ref{eq: sys}).

\demo{{\rm 4.1.} Quasiconvex functions} We begin by describing a quasi-convex function which will play
an important role in our construction using notation
introduced in Section~2.
We define $\fo\colon\sdd\to\RR$ by $\fo(X)=\det X$ when $X$
is positive definite and by $\fo(X)=0$ otherwise.
\spn{4.1}\proclaim{Lemma}\label{fo}
The function $\fo$ is quasiconvex on $\sdd$.
\endproclaim
 
{\it Proof}. This result is proved in \cite{Sv 92b}. In that paper the proof
is actually carried out for a more general class of functions.
We give a simple version of the proof here, for the convenience
of the reader. Let $\Omega=\{x\in\RR^2,\,|x|<1\}$ and let
$\phi\colon\Omega\to\RR$ be smooth and compactly supported in $\Omega$.
 We must prove that for each $A\in\sdd$ we have
$\int_{\Omega}(\fo(A+\nabla ^2\phi)-\fo(A))\ge0$.
This is obvious if $A$ is not positive definite, since then we integrate
a nonnegative function.
If $A$ is positive definite, we can assume $A=I$ by a simple change of 
variables. Let $u_0(x)=|x|^2/2$ and $u(x)=u_0(x)+\phi(x)$.
We also set $\varphi=\nabla u$, which will be viewed as a map
$\varphi\colon\Omega\to\RR^2$. Finally, we let 
$E=\{x\in\Omega,\, \det\nabla\varphi(x)\ge0\}$.
We must prove that $\int_{E}\det\nabla\varphi\ge\meas(\Omega)$.
Since $\det\varphi\ge0$ on $E$, we can use the area formula (\cite{Fe 69})
to infer that it is enough to prove $\Omega\subset\varphi(E)$.
Consider an arbitrary $b\in\Omega$ and let $a\in\bar\Omega$ be a point where
the function $x\to u(x)-b\cdot x$ attains its minimum in $\bar\Omega$.
It is easy to verify that $a\in \Omega$ and hence
$\varphi(a)=b$ and $a\in E$.
We see that $\Omega\subset\varphi(E)$ and the proof is finished.
\vglue8pt

In what follows we will use the following notation: for 
$X\in\mdd$ we let $X_\sym=(X+X^t)/2$ and $X_\as=(X-X^t)/2$.
\spn{4.2}\proclaim{Lemma}\label{extsym}
Let $f\colon\sdd\to\RR$ be a smooth function such that
$|D^2f|\le c$  in $\sdd$.
Assume that $f$ is strongly quasi\/{\rm -}\/convex in the sense that for some $\gamma>0$
we have
$\int_{\RR^2}(f(A+\nabla ^2\phi)-f(A))\ge\gamma\int_{\RR^2}|\nabla^2\phi|^2$
for all smooth{\rm ,} compactly supported $\phi\colon\RR^2\to\RR$.
Then for sufficiently large $\kappa>0$
the function $\tf\colon\mdd\to\RR$  defined
by $\tf(X)=f(X_\sym)+\kappa|X_\as|^2$ is strongly quasi\/{\rm -}\/convex.
\endproclaim
 \vglue-6pt
{\it Proof}.
Let $\Td$ be the two-dimensional torus $\RR^2/{\mathbf Z}^2$.
Let $\vf\colon\Td\to\RR^2$ be a smooth function and let
$A\in\mdd$. We want to prove that
$$\itd(\tf(A+\nabla \vf)-\tf(A))\ge \gamma/2\itd|\nabla \varphi|^2.$$
Let us consider the Helmholtz decomposition 
$\vf=\nabla \phi+\dper\eta+a$ of $\vf$, where $\phi$ and $\eta$
are scalar functions, $\dper\eta=J\nabla \eta$ (with $J$ as above),
and $a$ is a constant vector. 
We have $\nabla \vf=\nabla ^2\phi+\nabla \dper\eta$. 
Set $Y=(\nabla \dper\eta)_\sym$.
A standard  calculation (involving integration by parts and the use
of the identity $\itd\det \nabla ^2\eta\break =0$) gives
 $\itd|Y|^2=\itd  |\nabla
^2\eta|^2/2=\itd(\Delta\eta)^2/2=\itd|(\nabla\dper\eta)_\as|^2$.$\phantom{1^{\textstyle |^|}}$\hskip-6pt We can write
\begin{eqnarray*}
&&\hskip-64pt \itd(\tf(A+\nabla \vf)-\tf(A))   \\
&   =& \itd (f(A_\sym + \nabla ^2\phi)-f(A_\sym))  \\
& &  + \itd(\kappa|A_\as+(\nabla \dper\eta)_\as|^2-\kappa|A_\as|^2)  \\
& &  + \itd(f(A_\sym+\nabla ^2\phi+Y) - f(A_\sym + \nabla ^2\phi))  \\
&   =& I+II+III.
\end{eqnarray*}

We have $I\ge\gamma\itd|\nabla ^2\varphi|^2$ 
by our assumptions and Lemma~\ref{per}. The second term can be evaluated as
$II=\itd \kappa|Y|^2$ by using the calculation above and the fact that
$\itd \nabla ^2\eta=0$. Finally, the third term can be written as
\begin{eqnarray*}
III & =  & \itd(f(A_\sym+\nabla ^2\phi+Y)
 - f(A_\sym + \nabla ^2\phi)-Df(A_\sym+\nabla^2\phi)Y)  \\
& & +   \itd(Df(A_\sym+\nabla ^2\phi)-Df(A_\sym))Y  \\
&  \ge & -\itd(c/2|Y|^2+c|\nabla ^2\phi||Y|)   \\
& \ge &  -\itd(\gamma/2|\nabla ^2\phi|^2 + c/2|Y|^2 + c^2/(2\gamma)|Y|^2).
\end{eqnarray*}
We obtain the right inequality when 
$\kappa\ge \gamma/2+c/2+c^2/(2\gamma)$. The proof is finished.
\vglue6pt

Lemma~\ref{extsym} cannot be directly applied to the function $\fo$
from Lemma~\ref{fo}. However, we can modify $\fo$ in the following way.
We consider a smooth mollifier $\omega$ on $\sdd$ which is supported
in the  ball of radius $1/8$ centered at $0$ and satisfying 
$\ints \om = 1$, $\ints X\om(X)\,dX=0$, and $\ints\det(X)\om(X)\,dX=0$.
We let $f_1(X)=\max(\fo(X), |X|^2-25)$  and
$f_2=f_1*\om$.
We note that $f_2(X)=f_0(X)$ 
when $|X|\le 5$ and the open ball $B_{X,\frac18}$ is contained in the set of the
positive definite matrices.
Choosing a small $\gamma>0$ (to be specified later) and setting $f_3(X)=f_2(X)+\gamma|X|^2$, we denote by $\tf_3$
the strongly  quasi-convex extension of $f_3$ to $\mdd$ obtained in Lemma~\ref{extsym} (for a suitable $\kappa$).

Let $T=\left(\begin{array}{rr}
         0 & 1\\
         1 & 0
         \end{array}\right)$.
We define $\theta\colon\mdd\to\mdd$ by $\theta\cdot X=TXJ^t$, where $J$ is the rotation
by $\pi/2$ introduced above. Note that the diagonal matrices are invariant under
$\theta$ and that $\theta$ restricted to the diagonal matrices can be thought of as
 a rotation by $\pi/2$. The same is true for anti-diagonal matrices, by which we mean
the matrices of the form $TX$, where $X$ is diagonal. Therefore $\theta^2=-{\rm   Id}$.

We now define a function $f_4\colon\mdd\to\RR$, which will play an 
important r\^ole in our construction.
Let $H=\left(\begin{array}{rr}
         \frac54 & 0\\
         0 & -\frac54
         \end{array}\right)$,
and set
$$f_4(X)=\sum_{k=0}^3 \tilde f_3(\theta^{-k}\cdot X-H).$$
It is easy to see that $f_4$ satisfies $f_4(\theta\cdot X)=f_4(X)$ 
for each $X\in\mdd$ and therefore
$Df_4(\theta\cdot X)= \theta\cdot Df_4(X)$ for each $X\in\mdd$.
(We note that the restriction of $f_4$ to the diagonal matrices 
 vanishes in the square given by the matrices
$\theta^k\cdot H$, $k=0,1,2,3$, and on the
half-lines originating at $\theta^k\cdot H$ and passing through
$\theta^{k+1}\cdot H$, where $k=0,1,2,3$.)

We now let 
$$
A_1=\left(\begin{array}{rr}
         3 & 0\\
         0 & -1
         \end{array}\right),\,
A_2=\left(\begin{array}{rr}
         1 & 0\\
         0 & 3
         \end{array}\right),\,
A_3=\left(\!\begin{array}{rr}
         -3 & 0\\
         0 & 1
         \end{array}\right),\,
A_4=\left(\!\begin{array}{rr}
         -1& 0\\
         0 & -3
         \end{array}\right)\,,
$$
  noting that $A_{k+1}=\theta^k\cdot A_1$, $k=1,2,3$.
By a direct calculation,
$Df_4(A_1)=\left(\begin{array}{rr}
         \frac14 + 14\gamma & 0\\
         0 & \frac74 + 2\gamma
         \end{array}\right)$. 
By considering functions of the form 
$\frac12\alpha|X|^2+\beta f_4(X)$
 we can easily obtain the following lemma, by choosing suitable positive $\alpha, \beta$, and~$\gamma$. 

\spn{4.3}\proclaim{Lemma}\label{fff} \hskip-8pt
There exist a smooth{\rm ,} strongly quasi\/{\rm -}\/convex function $F_1\colon\mdd\break\to\RR$
with uniformly bounded $D^2F_1$ which satisfies {\rm (}\/in the notation introduced
above\/{\rm )}  $F_1(\theta\cdot X)=F_1(X)$ for each $X$ and
 $DF_1(A_1)=\left(\begin{array}{rr}
         1 & 0\\
         0 & 3
         \end{array}\right)$.
\endproclaim

{\it Proof}. See above.

\vglue12pt 
The set $K$ corresponding to the function $F=F_1$ (see the beginning
of the section) contains the matrices 
$ \left(\!\!\begin{array}{c}
         A_k \\
         DF_1(A_k)J
         \end{array}\!\!\right),\, k=1,\dots,4.$ 
These are the matrices
\begin{eqnarray*}
M_1^0=\!\left(\!\begin{array}{rr}
         3 & 0\\
         0 & -1\\
         0 & -1\\
         3 & 0
         \end{array}\!\right)\!,\, 
M_2^0=\!\left(\!\begin{array}{rr}
         1 & 0\\
         0 & 3\\
         0 & 3\\
         1 & 0
         \end{array}\!\right)\!,\, 
M_3^0=\!\left(\!\!\begin{array}{rr}
         -3 & 0\\
         0 & 1\\
         0 & 1\\
         -3 & 0
         \end{array}\!\right)\!,\,
M_4^0=\!\left(\!\!\begin{array}{rr}
         -1& 0\\
         0 & -3\\
         0 & -3\\
         -1 & 0
         \end{array}\!\right).\\
\end{eqnarray*}

\demo{{\rm 4.2.} Deformations of $T_4$-configurations} Let us consider four $m\times n$ matrices $M_1,\dots,M_4$.
We say that
$M_1,\dots,M_4$ are in $T_4$-configuration (see Figure 1) if
$\rank(M_i-M_j)\ne1$ for all $i,j$, and if there exist
rank-one matrices $C_1,\dots,C_4$ with $\sum_k C_k=0$,
real numbers $\kappa_1,\cdots\kappa_4>1$, and a \pagebreak matrix $P\in\mmn$
such that
\begin{eqnarray*}
M_1 & = & P+\kappa_1 C_1\,,\\
M_2 & = & P+C_1+\kappa_2 C_2\,,\\
M_3 & = & P+C_1+C_2+\kappa_3 C_3\,,\\
M_4 & = & P+C_1+C_2+C_3+\kappa_4 C_4.
\end{eqnarray*}

This configuration was discovered independently by several authors.
We are aware of  \cite{Sch 74}, where it is used in a similar context
as below,  \cite{AH 86}, and \cite{Ta 93}, where it is used in a 
different context.  
Slightly different examples 
exhibiting similar features were also independently discovered
 in \cite{NM 91} and \cite{CT 93}.  The paper \cite{BFJK 94} contains
an interesting example using a $T_4$-configuration. 
The following observation
appears in \cite{AH 86}, \cite{Ta 93} and implicitly 
also in the other papers.

\figin{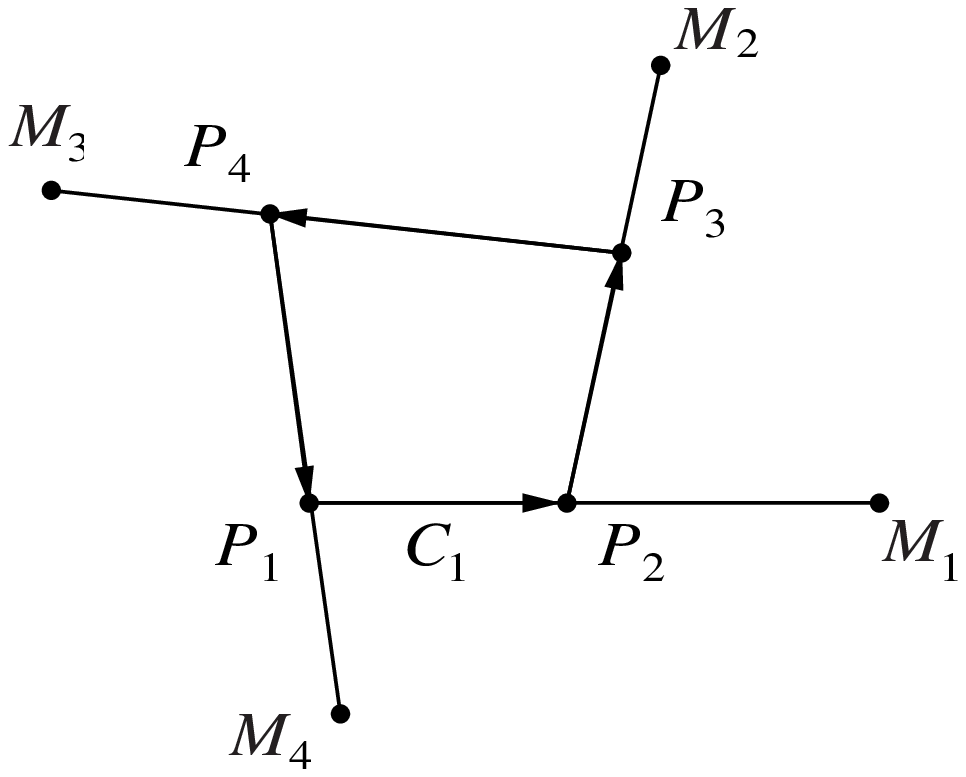}{570}
\begin{quote} Figure 1. A $T_4$ configuration with $P_1 = P$, $P_2 = P + C_1$,
$P_3 = P + C_1 + C_2$, $P_4 = P + C_1 + C_2 + C_3$. The lines
indicate rank-$1$ connections. Note that the figure need not be planar.
\end{quote}

\spn{4.4}\proclaim{Lemma}\label{tartar} 
If $M_1,\dots,M_4$ are in $T_4$\/{\rm -}\/configuration{\rm ,} the rank\/{\rm -}\/one convex
hull of the set $\{M_1,\dots,M_4\}$ contains the points
$P_1=P, P_2=P+C_1,P_3= P+C_1+C_2, P_4=P+C_1+C_2+C_3$.
For each point $X$ in the rank\/{\rm -}\/one convex hull there
exists a unique laminate $\mu = \sum \mu_l \delta_{M_l}$ with
center of mass $X$. 
\endproclaim

\demo{{P}roof}
To see this, let us consider a rank-one convex function $f\colon\mmn\break\to\RR$
which vanishes at the points $M_1,\dots,M_4$.
We have $$f(P_{i+1})\le 1/\kappa_i f(M_i)+(1-1/\kappa_i)f(P_i)=(1-1/\kappa_i)f(P_i)$$ for each
$i$, where the indices are considered modulo $4$. Applying this recursively, we get
that $f(P_i)\le 0$ for each $i$. 
Uniqueness is obvious if the $M_l$ span a three dimensional affine space. 
If all four matrices lie in a plane one can introduce coordinates
$x,y$ along the rank-one directions in this plane and exploit the fact that
the function $g(x,y)=xy$ satisfies $\langle\mu,g\rangle=g(\bar\mu)$.
\enddemo

{\it Example}. For   future reference, let us calculate the coefficients 
$\mu_l$ above for
$X=P_1$. We let $\beta_i=1-1/\kappa_i$, $i=1,\dots,4$.
Using recursively the identity $P_{i+1}=(1-\beta_i)M_i
+\beta_i P_i$ (where the indices are considered modulo $4$),
we get easily the following expression for the laminate $\mu$ supported
in $\{M_1,\dots,M_4\}$ with $\bar\mu=P_1$:
\begin{equation}\label{laminate}
\mu=\sum_{i=1}^{4} 
\frac{(1-\beta_i)\beta_1\beta_2\beta_3\beta_4}
{\beta_1\dots\beta_i(1-\beta_1\beta_2\beta_3\beta_4)}\delta_{M_i}\,.
\end{equation}

\bigskip

The matrices $M_k^0$ at the end of subsection~4.1 
 are in $T_4$-configuration, as one can see by taking
\begin{eqnarray*}
P=\!\left(\begin{array}{rr}
         -1 & 0\\
         0 & -1\\
         0 & -1\\
         -1 & 0
         \end{array}\right),\,\, 
C_1=\!\left(\begin{array}{rr}
         2 & 0\\
         0 & 0\\
         0 & 0\\
         2 & 0
         \end{array}\right), \,\,
C_2=\!\left(\!\begin{array}{rr}
         0 & 0\\
         0 & 2\\
         0 & 2\\
         0 & 0
         \end{array}\right),
\end{eqnarray*}
and $C_3=-C_1, \,C_4=-C_2, \,\kappa_1=\kappa_2=\kappa_3=\kappa_4=2$.
The matrices also lie in the set 
\begin{eqnarray*}
K_1=
\left\{\left(
\begin{array}{c} X \\ DF_1(X)J \end{array}\right);\, X\in\mdd\right\}\subset\mcd
\end{eqnarray*}
given by the quasi-convex 
function $F_1$ constructed
in Lemma~\ref{fff}. This shows that the rank-one convex hull $\Kjrc$ of $K_1$
is nontrivial. We now wish to establish that $\Kjrc$ is sufficiently large,
so that we can apply Theorem~\ref{main}.
We will see later that rather than trying to work with the specific function 
$F_1$, it is more convenient to work with a small perturbation
$F=F_1+\ve V$ of $F_1$, where $V$ is a compactly supported smooth function,
the properties of which will be specified later. For the moment we will only
assume that $F$ satisfies $DF(A_k)=DF_1(A_k)$ for $k=1,2,3,4$, where the matrices
$A_k$ are the same as in Subsection~4.1 .
We also denote by $K\subset\mcd$ the set corresponding to $F$.
By our assumptions we know that $K$ contains a $T_4$-configuration
given by the matrices $M^0_k$, $k=1,2,3,4$ defined above.
It is natural  to investigate deformations of this
$T_4$-configuration.  In other words, we will investigate four-tuples 
$M_1,\dots M_4$ such that, for $k=1,\dots,4$,
$M_k$ is close to $M_k^0$, $M_k\in K$, and $M_1,\dots M_4$ are in $T_4$-configuration.
\vglue12pt
\def\nhs{\hskip-4pt}
We introduce the following notation.
$$
  \begin{array}{rlrll}
e_1  &\nhs= (1,0)\,,& \qquad 
e_2 &\nhs =  (0,1)\,, \\
f_1 &\nhs=  (2,0,0,2)\,, & 
f_2 &\nhs=  (0,2,2,0)\,,\\
C_1^0 &\nhs =  f_1\otimes e_1\,,  &
C_2^0  &\nhs=  f_2\otimes e_2\,,\\
C_3^0& \nhs=  -C_1^0\,,  &
C_4^0  &\nhs=  -C_2^0\,,\\
P^0  &\nhs=  -(C_1^0+C_2^0)/2\,,   \\
 \kappa_1^0&\nhs  =  \kappa_2^0 = 
\kappa_3^0 = \kappa_4^0  = 2\,.
\end{array}
$$  

We parametrize the rank-one matrices $C_k$ in a small neighborhood of 
$C_k^0$ as follows.
\begin{eqnarray*}
C_1 & = & (f_1+a_1)\otimes(e_1+\beta_1e_2)\,,\\
C_2 & = & (f_2+a_2)\otimes(e_2-\beta_2e_1)\,,\\
C_3 & = & (-f_1+a_3)\otimes(e_1+\beta_3e_2)\,,\\
C_4 & = & (-f_2+a_4)\otimes(e_2-\beta_4e_1)\,,
\end{eqnarray*}
where $a_1,\dots ,a_4$ are (small) vectors in $\RR^4$, 
and $\beta_1,\dots, \beta_4$ are (small) real numbers.
We linearize the equation $\sum_k C_k=0$ around the
solution $C_k^0$. The linearized equation is equivalent to
\begin{eqnarray*}
a_1+a_3 + (\beta_4-\beta_2)f_2 & = & 0\,,\\
a_2+a_4 + (\beta_1-\beta_3)f_1 & = & 0\,.
\end{eqnarray*}
Using these formulae and the above expressions for $M_k$, we easily check
(with the help of the implicit-function theorem) that the four-tuples
$(M_1,\dots,M_4)$ of the $4\times 2$  matrices which are close to 
$(M_1^0,\dots,M_4^0)$ and form $T_4$-configuration such that
the parameters $P,C_j,\kappa_j$
are close to $P^0,C_j^0,\kappa_j^0$
form a 24-dimensional manifold $\MM$. The tangent space $L_{\MM}$
of $\MM$ at the point $(M_1^0,\dots,M_4^0)$ can be identified with
 four-tuples $(Z_1,\dots,Z_4)$ of $4\times2$ matrices
of the form 
\begin{eqnarray*}
Z_1 & = & \left(\begin{array}{ll}
         p_{11}+2a_{11}+\varkappa_1' & p_{12}+2\beta_1'\\
         p_{21} +2a_{21} & p_{22}\\
         p_{31}+2a_{31} & p_{32}\\
         p_{41}+2a_{41}+\varkappa_1'& p_{42}+2\beta_1'
         \end{array}\right), \\[4pt]
Z_2 & = & \left(\begin{array}{ll}
         p_{11}+a_{11} & p_{12}+2a_{12}+\beta_1'\\
         p_{21} +a_{21}-2\beta_2' & p_{22}+2a_{22}+\varkappa_2'\\
         p_{31}+a_{31}-2\beta_2' & p_{32}+2a_{32}+\varkappa_2'\\
         p_{41}+a_{41}& p_{42}+2a_{42}+\beta_1'
         \end{array}\right), \\[4pt]
Z_3 & = & \left(\begin{array}{ll}
         p_{11}-a_{11}-\varkappa_3' & p_{12}+a_{12}-2\beta_3'+\beta_1'\\
         p_{21}-a_{21}+\beta_2'-2\beta_4' & p_{22}+a_{22}\\
         p_{31}-a_{31}+\beta_2'-2\beta_4' & p_{32}+a_{32}\\
         p_{41}-a_{41}-\varkappa_3'& p_{42}+a_{42}-2\beta_3'+\beta_1'
         \end{array}\right), \\[4pt]
Z_4 & = & \left(\begin{array}{ll}
         p_{11} & p_{12}-a_{12}+\beta_3'-\beta_1'\\
         p_{21} +\beta_4' & p_{22}-a_{22}-\varkappa_4'\\
         p_{31}+\beta_4' & p_{32}-a_{32}-\varkappa_4'\\
         p_{41} & p_{42}-a_{42}+\beta_3'-\beta_1'
         \end{array}\right),
\end{eqnarray*}
where the values of all the 24 parameters run through all real numbers. 
Moreover, there \pagebreak is a well-defined mapping 
$(M_1,\dots,M_4)\to (P_1,\dots,P_4)$ from $\MM$ to the
four-tuples of $4\times 2$ matrices, where (in the notation introduced
in the definition of $T_4$-configuration)
$P_1=P,\,P_2=P_1+C_1,\,P_3=P_2+C_2,\,P_4=P_3+C_3$ as above.

We now consider the additional constraint $M_k\in K$, where $K$ is the
set determined by $F$. The four-tuples $(M_1,\dots,M_4)$ 
satisfying $M_k\in K$ clearly form a 16-dimensional manifold 
$\KK=K\times K\times K\times K$. The tangent space $L_{\KK}$
of $\KK$ at $(M_1^0,\dots,M_4^0)$ can be identified with the four-tuples
\begin{eqnarray*}
\left(\!\begin{array}{c} X_1 \\ D^2F(A_1)X_1J \end{array}\!\right),
\left(\!\begin{array}{c} X_2 \\ D^2F(A_2)X_2J \end{array}\!\right),
\left(\!\begin{array}{c} X_3 \\ D^2F(A_3)X_3J \end{array}\!\right),
\left(\!\begin{array}{c} X_4 \\ D^2F(A_4)X_4J \end{array}\!\right)
\end{eqnarray*}
where $X_1,\dots,X_4$ run through all $2\times 2$ matrices.

We now 
consider the maps   $(M_1,\dots,M_4)\to (M_k,P_k')$, where 
 $P_k$ is defined as
above and where we denote (with a slight
abuse of notation)  by $P_k'$ the orthogonal projection
of the point $P_k$ into the space $(T_{A_k} K)^\perp$,
the normal space of $K$ at $A_k$.
We would like to establish  the following nondegeneracy
conditions, which will be important later when we construct
in-approximations.

\demo{Condition {\rm (C)}}  $\MM$ and $\KK$ intersect
transversely at $(M_1^0,\dots,M^0_4)$ and,\break
 (after $\MM$ is perhaps replaced by a sufficiently small neighborhood of\break
$(M_1^0,\dots,M_4^0)$  in $\MM$) the map
 $(M_1,\dots,M_4)\to (M_k,P_k')$ is, for each $k$,  a nondegenerate
diffeomorphism of  $\MM\cap\KK$ and a neighborhood of 
$(M_k^0, (P_k^0)')$
in $K\times(T_{A_k}K)^\perp$.\enddemo

Rather than trying to decide whether these nondegeneracy conditions are
satisfied for an explicitly given function $F$, it seems to be more natural
to verify that the conditions are satisfied in the generic case.
More specifically, we note that for each smooth compactly
supported function $V\colon\mcd\to\RR$ the function $F=F_1+\varepsilon V$
is strongly quasi-convex for sufficiently small $\ve$.
By choosing $V$ in a suitable way, we can perturb $D^2F(A_1),\dots
D^2F(A_4)$ to any prescribed values which are close enough
to the original values, without changing the values of $DF(A_1),\dots,
DF(A_4)$, and without affecting the strong quasi-convexity. 
For the purpose of the construction of the counterexample 
announced at the beginning
of this section, we can therefore restrict our considerations to the
generic case.

\spn{4.5}\proclaim{Lemma}\label{tranversal}
Assume that $DF(A_k)=DF_1(A_k)$ for $k=1,2,3,4$.
Then condition {\rm (C)} above is satisfied for the 
generic values of $D^2 F(A_k)${\rm ,} $k=1,\dots,4$.
\endproclaim
\vglue-14pt

{\it Proof}.
The condition that $\MM$ and $\KK$ intersect
transversely at $(M_1^0,\dots,M^0_4)$ and that
  the map
 $(M_1,\dots,M_4)\to (M_1,P_1')$ is a nondegenerate
diffeomorphism of a small neighborhood of $(M_1^0,\dots,M^0_4)$ in
$\MM\cap\KK$ and a neighborhood of 
$(M_1^0, (P_1^0)')$
in $K\times(T_{A_1}K)^\perp$ is easily seen to be equivalent to the condition that
the following linear homogeneous system of 40 equations for 40 unknowns has no
nontrivial \pagebreak solutions.
\begin{eqnarray*}
 Z_j & = & \left(\begin{array}{c} X_j \\ D^2F(A_j)X_jJ \end{array}\right), \quad j=1,2,3,4,\\
\left(\begin{array}{ll}
         p_{31} & p_{32}\\
         p_{41} & p_{42}\\
         \end{array}\right) &= &
 D^2F(A_1)\left(\begin{array}{ll}
         p_{11} & p_{12}\\
         p_{21} & p_{22}\\
         \end{array}\right)J,\\
X_1 &= &0,
\end{eqnarray*}
where $Z_j=Z_j(p_{kl},a_{kl},\beta_k',\kappa_k')$ (with $k=1,2,3,4,\,l=1,2$)
 are the $4\times 2$ matrices
introduced above and $X_1,X_2,X_3,X_4$ are $2\times 2$ matrices.
The determinant of the corresponding $40\times 40$ matrix 
is a polynomial expression in the entries of the matrices $D^2 F(A_j)$
(which are now considered as parameters), and will be denoted by $Q_1$. 
The polynomial $Q_1$ is not identically zero, since for
 $$D^2F(A_1)=I,\quad D^2F(A_2)=I,\quad D^2F(A_3)=0,\quad D^2F(A_4)=I$$
 we can
 check by a straightforward calculation
that the system has no nontrivial solutions.

By using symmetry we see that, for each $k=1,2,3,4$,
the condition that $\MM$ and $\KK$ intersect
transversely at $(M_1^0,\dots,M^0_4)$ and that
  the map
 $(M_1,\dots,M_4)\break\to (M_k,P_k')$ is a nondegenerate
diffeomorphism of a small neighborhood of\linebreak $(M_1^0,\dots,M^0_4)$ in
$\MM\cap\KK$ and a neighborhood of 
$(M_k^0, (P_k^0)')$
in $K\times(T_{A_k}K)^\perp$ can be expressed as $Q_k\ne 0$, where
$Q_k$ is a suitable nonzero polynomial 
in the entries of the matrices $D^2 F(A_j)$.
Hence all of our nondegeneracy conditions will be satisfied
at all values of $D^2F(A_j)$ where the polynomial $Q=Q_1Q_2Q_3Q_4$ does not vanish.
Since $Q$ is not identically zero, the result follows.

\demo{{\rm 4.3.} In-approximation}
To be able to use Theorem~\ref{main}, we need to have a suitable
in-approximation.  \enddemo

\spn{4.6}\proclaim{Lemma}\label{inapp}
Using the notation above{\rm ,}  
assume that condition {\rm (C)} is satisfied. Let $r>0$. Then 
there exists an in\/{\rm -}\/approximation
 $\{\UU_i\}_{i=1}^\infty$ of 
 $$K_r=\mathbold{\cup}_{j=1}^4 \{X\in M^{4\times 2},\, |X-M^0_j|\le r\}\cap K$$ 
 such that $\UU_1$ contains a {\rm (}\/small\/{\rm )} neighborhood of 
the rank\/{\rm -}\/one convex hull of the points $P^0_1,\dots,P^0_4$. 
\endproclaim

{\it Proof.} 
 Let $\OO$ be a sufficiently small neighborhood of
$(M^0_1, M^0_2, M^0_3, M^0_4)$ in 
$\MM \cap \KK \subset (\mcd)^4$.
The main point is that, for each $k=1,2,3,4$, the image of $\OO$ under
\newcommand{\MMM}{{(M_1,M_2,M_3,M_4)}}
 the map $\MMM\to P_k\MMM$ is open in $\mcd$,
 whereas the image of $\OO$ under the projections $\MMM\break\to M_k$
  is not
(since $M_k\in K$). We will therefore
consider convex combinations
$(1 - \lambda) P_k +  \lambda  M_k$ with $\lambda \to 1$
to construct an in-approximation of $K$ (see Fig.~2).
We now describe the details.

\centerline{\BoxedEPSF{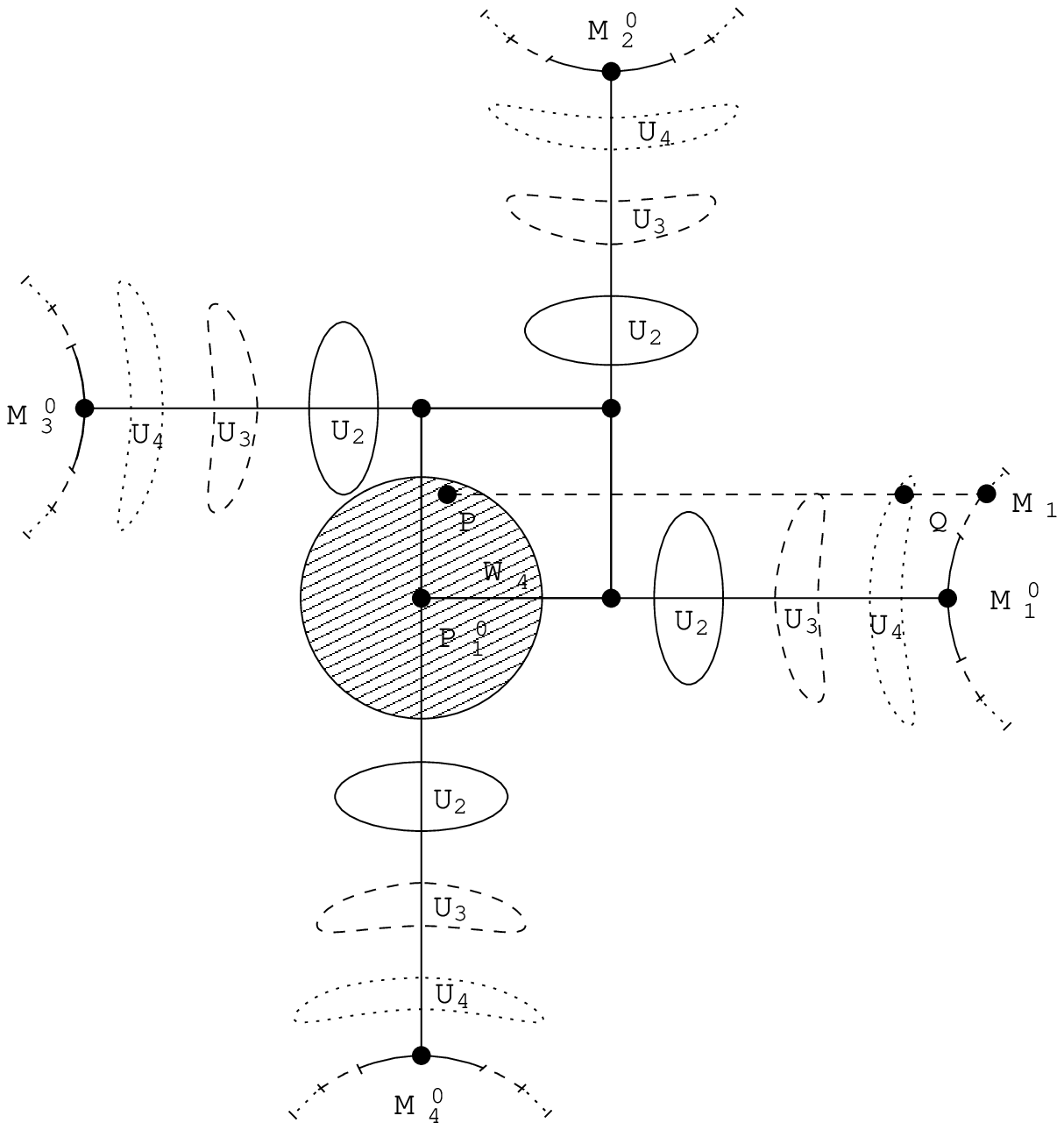 scaled 750}}
\begin{quote} Figure 2. Schematic illustration of the sets $\UU_2, \UU_3, \UU_4 \subset
\mcd$. The  solid (resp.\ dashed, or dotted) lines
through the point $M^0_1$ are the projections of the set
 $\OO_2 \mbox{(resp.\ $\OO_3$, or $\OO_4$)} \subset \MM\cap\KK\subset(\mcd)^4$
to the first component. They are not open in $\mcd$ since
they are contained in $K$. The shaded set $W_4$ is 
the image of $\OO_4$
under the map $(M_1,M_2,M_3,M_4)\to P_1(M_1,M_2,M_3,M_4)$ 
and it is open in $\mcd$. 
By $P = P_1(M_1,M_2,M_3,M_4)$ we  denote a typical point in $W_4$.
A typical point  $Q$  in $\UU_{1,4} \subset \UU_4$
is given by\break $(1 - \lambda_4) P_1(M_1,M_2,M_3,M_4)+  
\lambda_4 M_1$, where $(M_1,M_2,M_3,M_4) \in \OO_4$.
\end{quote}

 We consider a sequence 
$\OO_0,\OO_1,\OO_2\dots\subset\OO$
of open neighborhoods of \linebreak $(M^0_1,\dots,M^0_4)$ in $\MM\cap\KK$,
such that each $\OO_j$ is diffeomorphic to the eight-dimensional
unit ball and that, for each $j=0,1,2,\dots$ we have
$\bar\OO_j\subset\OO_{j+1}$. We  also consider a sequence of numbers
$0=\lambda_0$,  $1/2<\lambda_1<\dots<\lambda_j<\dots<1$
converging to $1$
as $j\to\infty$. For $j=0,1,2,\dots$ we let
$$\UU_{k,j}=\{(1-\lambda_j)P_k+\lambda_j M_k,\,
(M_1,\dots,M_4)\in\OO_j\},$$ where $P_k=P_k(M_1,\dots,M_4)$ 
is the map considered in subsection~4.2.
We also let $\UU_j=\cup_{k=1}^{k=4} \UU_{k,j}$.
Condition (C) implies  that there exists
$j_0$ such that the sets $\UU_j$
are open when $j\ge j_0$ and $\OO$ is sufficiently small.
To see this, consider for example $k=1$ and let us write 
points $M_1\in K$ which are close to $M^0_1$ as
$M_1=M_1^0 + X+\xi(X)$, with $X\in T_{A_1}K$ and $\xi(X)\in (T_{A_1}K)^\perp $.
We can also write $P_1=P^0_1+Y+\eta$ with $Y\in (T_{A_1}K)^\perp$ and
$\eta\in T_{A_1}K$.
If condition (C) is satisfied, we know that, in a small neighborhood
of $(M_1^0,\dots,M_4^0)$, we can take $X$ and $Y$ as local
coordinates in $\MM \cap \KK$. For $(M_1,\dots,M_4)\in\MM\ \cap \KK$ which is
close to $(M_1^0,\dots,M_4^0)$ and $P_1=P_1(M_1,\dots,M_4)$,
we can therefore write the $\eta-$component of $P_1$ in the above
decomposition as
$\eta=\eta(X,Y)$, where $\eta$ is a smooth function of 
$X$ and $Y$ with $\eta(0,0)=0$.
In the coordinates $(X,Y)$, the derivative of the map
$(X,Y)\to(1-\lambda)P_1+\lambda M_1$ is given
by  the block matrix
$$
\left(\begin{array}{cc}
         \lambda I+(1-\lambda) \partial_X\eta &  (1-\lambda)\partial_Y\eta \\
         \lambda\partial_X\xi & (1-\lambda)I
         \end{array}\right)\,.
$$
Since $\partial_X\xi(0)=0$, we see that the matrix is regular 
when $X$ is small and $\lambda$ is close to $1$.
The openess of $\UU_{1,j}$ for large $j$, $\lambda$ close (but not equal)
 to $1$, and  small $\OO$
follows.

By Lemma~\ref{splitting} below, the closure of $\UU_j$ (and hence the
closure of its rank-one convex hull) is contained in the rank-one convex
hull of $\UU_{j+1}$.  
 Moreover, the rank-one
convex hull of $\UU_0$ contains a neighborhood of the
square given 
by the convex hull of the points $P^0_1,\dots,P^0_4$
(which coincides with the rank-one convex hull of these points,
since the points lie in a two-dimensional plane).
The required in-approximation has therefore been established.

\spn{4.7}\proclaim{Lemma}\label{splitting} Using the notation introduced in the proof
of Lemma~{\rm \ref{inapp}} the following is true. For each $j=1,2,\dots\,${\rm ,}
the set $\overline{\UU_j}$ is contained in $\UU_{j+1}^{\mathrm rc}${\rm ,} and
each $A \in \overline{\UU_{j,k}}$ is the center
of mass of a laminate $\mu = \sum_{l=1}^4 \mu_l \delta_{Y_l}${\rm ,}
with $Y_l \in \UU_{l,j+1}$. Moreover{\rm ,} when  $\lambda_j$ is sufficiently close to $1$ and
$\OO$  is sufficiently small{\rm ,} we can
achieve in addition that
\begin{eqnarray*}
\mu_k &\geq& 1 - (\lambda_{j+1} - \lambda_j)\,,\\
 |Y_k - A| &\leq&
 2|M^0_1-P^0_1|(\lambda_{j+1} - \lambda_j)\,, \\
\mu_l &\geq & (\lambda_{j+1} - \lambda_j)/8, \quad \mbox{for} \, \, l \neq k. \\
\noalign{\vskip-28pt}
\end{eqnarray*}
\endproclaim

{\it Proof}.
 To simplify the notation suppose $A \in \overline{\UU_{1,j}}$. Then
there exist $(M_1,M_2,M_3,M_4) \in \overline{\OO_j} \subset 
\OO_{j+1}$ such that $A = (1- \lambda_j) P_1 + \lambda_j M_1$, where
$P_1=P_1\MMM$.
Let $Y_l = (1- \lambda_{j+1}) P_l + \lambda_{j+1} M_l$ (see Fig.\ 3).
Then $A$ is the center of mass of the laminate
\begin{eqnarray*}
\tilde{\mu} = \frac{\lambda_j}{\lambda_{j+1}} \delta_{Y_1} + 
   ( 1 -  \frac{\lambda_j}{\lambda_{j+1}}) \delta_{P_1}\,.
\end{eqnarray*}
and $|Y_1 - A|=|M_1-P_1| 
(\lambda_{j+1} - \lambda_j)\le 2|M^0_1-P^0_1|(\lambda_{j+1} - \lambda_j)$.

By Lemma~\ref{tartar} the point $P_1$ is the center
of mass of a unique laminate $\eta = \sum_{l=1}^4 \alpha_l \delta_{Y_l}$
supported on the $T_4$ configuration \pagebreak $(Y_1,Y_2,Y_3,Y_4)$, where\break

\vglue-32pt
\figin{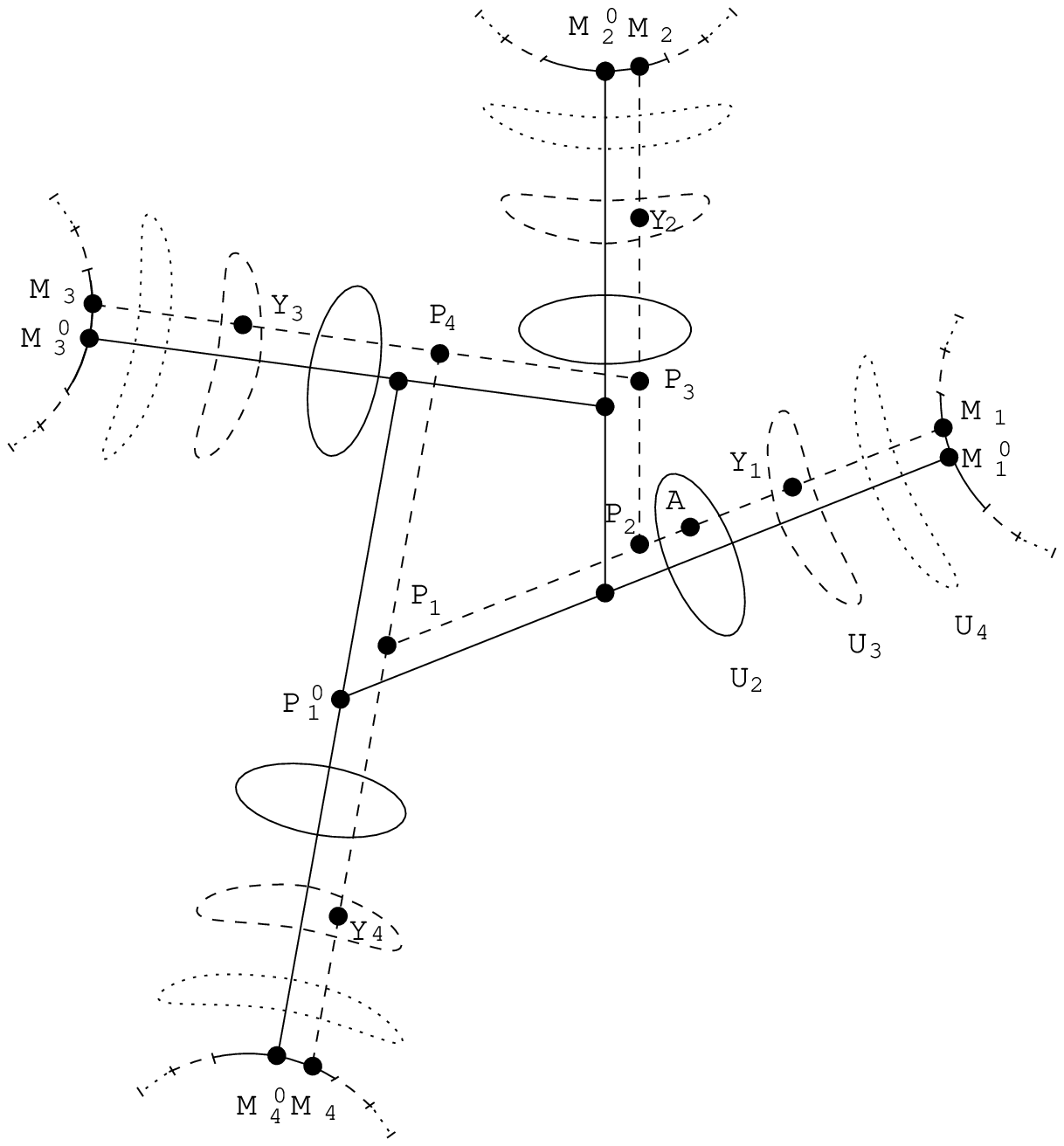}{750}
\begin{quote} Figure 3. The point $A \in \UU_{1,2}$ lies in the rank-1 convex hull
of the points $Y_1, Y_2, Y_3, Y_4 \in \UU_3$.
\end{quote}

\noindent the coefficients
$\alpha_l$ are given by  equation~(\ref{laminate}). 
Since $P_1$ and $Y_1$ differ by a rank-one matrix the measure
\begin{eqnarray*}
\mu = \frac{\lambda_j}{\lambda_{j+1}} \delta_{Y_1} +
( 1 -  \frac{\lambda_j}{\lambda_{j+1}}) \eta
\end{eqnarray*}
is a laminate with center of mass $A$. 
If $\OO$ is small and $\lambda_j$
is close to $1$, the numbers $\beta_i$ in (\ref{laminate}) are close to $1/2$,
and an elementary calculation gives our estimates.

\demo{{\rm 4.4.} Solutions with nowhere continuous gradients} 
\enddemo
{\it Proof of Theorem}~\ref{solutions}. 
The main idea of the proof is described in heuristic terms in the remarks
 immediately following the
theorem. In the proof below we will be freely using the notation introduced
earlier in Section 4. 

\smallskip
We recall that   $A_1,\dots, A_4$ are defined as follows:
$$
A_1=\left(\begin{array}{rr}
         3 & 0\\
         0 & -1
         \end{array}\right),\,
A_2=\left(\begin{array}{rr}
         1 & 0\\
         0 & 3
         \end{array}\right),\,
A_3=\left(\!\begin{array}{rr}
         -3 & 0\\
         0 & 1
         \end{array}\right),\,
A_4=\left(\!\begin{array}{rr}
         -1& 0\\
         0 & -3
         \end{array}\right)\,;
$$
see Section~4.1 .
We let $F_0$ be a suitable small perturbation of the quasiconvex function $F_1$ from
Lemma~\ref{fff} such that  $DF_0(A_k)=DF_1(A_k)$ for $k=1,\dots,4$ and
condition (C) is satisfied. Since the transversality and the other nondegeneracy conditions
are stable under small perturbations,
a version of (C) with $M^0_1,\dots,M^0_4$ replaced by 
close-by matrices $\tilde M^0_1\dots,\tilde M^0_4$ will also be satisfied
for any $F$ as in the statement of the theorem, provided $\delta$ is sufficiently small.
Moreover, we see easily that by choosing $\delta$ sufficiently small we can also achieve
that  Lemma~\ref{inapp} can be applied (with 
$M^0_1,\dots,M^0_4$ replaced by 
close-by matrices $\tilde M^0_1\dots,\tilde M^0_4$) with a fixed small $r>0$ 
to any set $K$ arising from a function $F$ satisfying
the assumptions of the theorem. In addition, we see easily that the in-approximations  
can be constructed 
so that $\UU_1$ contains a fixed small neighborhood of the zero matrix 
for any $F$ satisfying the assumptions. Let us choose $\varepsilon>0$
so that the ball of radius $\varepsilon$ centered at the zero matrix
is contained in this fixed small neighbourhood.
We see that the assumptions  of Theorem~\ref{main} are satisfied in our situation. 
However, it does not seem to be 
immediately clear that the solutions obtained from Theorem~\ref{main} are not
continuously differentiable in any open subset of $\Om$.
To obtain such solutions in a simple way, we  make the construction more
explicit and impose some additional conditions on the approximations
so that the nowhere differentiability of the limit is easy to see.

Let $\{\lambda_j\}$ and $r>0$ be as in Lemma~\ref{inapp}, 
and assume (as we can without loss of
generality) that $r$ is sufficiently small.
Let $\UU_j$ denote the in-approximation 
constructed in 
Lemma~\ref{inapp}.
Let $\phi\colon\mcd\to\RR$ be
be a continuous function which is $\equiv 1$ in $\{X;\, |X|\le 2r\}$ and vanishes
outside $\{X, |X|\le 3r\}$. For $l=1,2,3,4$ set $\phi_l(X)=\phi(X-M^0_l)$. 
Assume now that $\varepsilon$ is as above, $v\colon\Om\to\RR^2$ is as in
Theorem~\ref{solutions} and let $\varepsilon_1\colon\Om\to\RR$ be a continuous 
function in $\Omega$ which is $>0$. Let 
$\tilde w=\left(\begin{array}{l} v \\ 0 \end{array}\right)$.
We will now go through constructions involved in the
proof of Theorem~\ref{main} in more detail and construct a sequence of functions
$w_j\colon\Om\to\RR^4$ together with a sequence $\FF_j$ of families of open subsets of 
$\Om$,  so that  the following conditions are satisfied.

\begin{itemize}
\item[(i)]      The sets in $\FF_j$ are open, mutually disjoint,  contained
                in $\Omega$ together with their closures, and cover $\Om$
                up to a set of measure zero;
\item[(ii)]     Each set of $\FF_{j+l}$ is contained in a set
                of $\FF_j$ (where $j,l\ge 1$);
\item[(iii)]    $\sup\,\{{{{\rm   diam}}}\, V;\, V\in\FF_j\}\to 0$ 
                as $j\to\infty$;
\item[(iv)]  $\nabla w_j$ is constant on $V$ for each $V\in\FF_j$;
\item[(v)]      $\nabla w_j\in\UU_j$ a.e.\  in $\Om$;
\item[(vi)] $|w_1-\tilde w|<\varepsilon_1/2$ in $\Om$ and 
              $|w_{j+1}-w_j|\le 2^{-j-2}\varepsilon_1$ in $\Omega$, ($j=1,2,\dots $)\,.
\end{itemize}
In addition, the following conditions, which are crucial for the
desired behavior,  are satisfied when $j$ is sufficiently
large.

\begin{itemize}
\item[(vii)] ({\it $L^1$-convergence of $\nabla w_j$}) We have
 $\int_{\Om}|\nabla w_{j+1}-\nabla w_j|\le 
            L(\lambda_{j+1}-\lambda_j)\,\meas\,\Om\,\,$
             for a suitable constant $L$;

\item[(viii)] ({\it Persistence of oscillations}) For each $V\in\FF_j$ and each $l\in\{1,2,3,4\}$,
\begin{eqnarray}\label{p1}
\int_V\phi_l(\nabla w_{j+1}) &\ge &\frac18(\lambda_{j+1}-\lambda_j)\,\meas\, V \quad\mbox{and}\\
\label{p2}
\int_V\phi_l(\nabla w_{j+1}) & \ge & (1-(\lambda_{j+1}-\lambda_j))\int_V\phi_l(\nabla w_{j})\,.
\end{eqnarray}
\end{itemize}

Once the existence of $\{w_j\}$ and $\{\FF_j\}$ satisfying (i)--(viii) is established,
we can consider $w_\infty=\lim_{j\to\infty}w_j$. From (v)--(vii) we infer that 
\newcommand{\wi}{{w_{\infty}}}
$\wi$ is Lipschitz, with $\nabla\wi\in K$ a.e.\  in $\Om$. Moreover, using (ii), (vii),
and (viii) we see that, for each sufficiently large $j$ and $V\in\FF_j$,
\begin{eqnarray*}
  &&\hskip-36pt\int_V \phi_l (\nabla\wi) =
  \lim_{m\to\infty} \int_V \phi_l(\nabla w_m) \\
             && \ge  \lim_{m\to\infty}
(1-(\lambda_m-\lambda_{m-1}))
\dots (1-(\lambda_{j+2}-\lambda_{j+1}))
\int_V \phi_l(\nabla w_{j+1}) \\
&& \ge \frac{1}{16}\lambda_{j+1}(\lambda_{j+1}-\lambda_j)\,\meas V\,.
\end{eqnarray*}
This, together with (iii) implies that the essential oscillation of
$\nabla\wi$ over any open set is at least $\max_{1\le k<l\le 4}|M^0_k-M^0_l|/2$,
and therefore $\wi$ cannot be continuously differentiable in any open subset
of $\Om$.

To construct $\{w_j\}$ and $\{\FF_j\}$, we proceed by induction. The existence
of $w_1$ and $\FF_1$ satisfying (i)--(v) and the first
inequality of (vi) follows from Theorem~\ref{open}.
Assume that, for some $j\ge 1$ there exist $w_j$ and $\FF_j$ satisfying 
(i), (iv), and (v). Let $V\in\FF_j$ and assume that $\nabla w_j=A$
in $V$, with $A\in\UU_j$. Assume
that $A\in\UU_{1,j}$, for example. By Lemma~\ref{splitting}, the matrix $A$ is the center
of mass of a laminate $\mu=\sum_{l=1}^4\mu_l\delta_{Y_l}$, with
$Y_l\in\UU_{l,j+1}$. In addition, by the same lemma, if $j$ is sufficiently 
large,  
\begin{eqnarray}
\label{c1}
|Y_1-A| & \le & 2|M^0_1-P^0_1|(\lambda_{j+1}-\lambda_j)\,, \\
\label{c2}
\mu_1 &\ge & 1-(\lambda_{j+1}-\lambda_j)\quad \mbox{and}\\
\label{c3}
\mu_l &\ge & (\lambda_{j+1}-\lambda_j)/8 \quad\mbox{for each $l=1,2,3,4$}.
\end{eqnarray}

From the remark following the proof of Theorem~\ref{open} we see that there
exists a piecewise affine 
\newcommand{\wjv}{{w^V_{j+1}}}
$\wjv\colon V\to\RR^4$ such that $\nabla\wjv\in\UU_{j+1}$ a.e.\break 
in $V$, $|\wjv-w_j|\le \varepsilon_1 2^{-j-2}$ in $V$, $\wjv=w_j$ at the boundary
of $V$, and\break
$\meas\{x\in V;\,\nabla\wjv\in\UU_{l,j+1}\}=\mu_l\,\meas V$ for each 
$l=1,2,3,4$.
We choose a family $\FF_{j+1}^V$ of mutually disjoint open sets of radius $<1/(j+1)$
 which cover
$V$ up to a set of measure zero and $\nabla\wjv$ is constant on each of them.
We can now define $\FF_{j+1}=\cup_{V\in\FF_j}\FF_{j+1}^V$ and
$w_{j+1}\colon\Omega\to \RR^4$ by $w_{j+1}=\wjv$ in the closure of $V$ for each
$V\in\FF_j$.

From (\ref{c1}) and (\ref{c2}) we see that (vii) is satisfied with
  $$L=2|M^0_1-P^0_1|+2\max_{1\le k<l\le 4}|M^0_k-M^0_l|.$$
In addition, (\ref{c2}) and (\ref{c3}) imply that (viii) is satisfied.
The rest of the properties (i)-(viii) are immediate consequences of
our construction.

\demo{{R}emark}
The above construction is quite similar to the following simpler
example. Let us consider a sequence $0<\lambda_0<\lambda_1<\dots
\lambda_j<\dots<1$, with $\lim_{j\to\infty} \lambda_j=1$.
Let $X\subset L^{\infty}(0,1)$ be the space of all piecewise
constant functions. For a function $f\in X$ with $|f|\le\lambda_j$
we define $T_jf\in X$ in the following way. Let $(a,b)$ be a maximal
open interval on which $f$ is constant. Let $c=(a+b)/2$. We 
find $d\in (a,c)$ and $e\in (c,b)$ such that 
the function $g\colon(a,b)\to \RR$ defined by
$g(x)=-\lambda_j$ when $x\in(a,d)$, $g(x)=\lambda_j$ when $x\in (d,c)$,
$g(x)=-\lambda_j$ when $x\in(c,e)$, and $g(x)=\lambda_j$ when $x\in(e,b)$
has the same average as $f$ over the intervals $(a,c)$ and $(c,b)$.
We then set $T_jf(x)=g(x)$ for $x\in(a,b)$, and repeat the same
construction on the other maximal intervals on which $f$
is constant. Let $0<A<\lambda_0$ and let $f_0\equiv A$ in $(0,1)$.
Set $f_{j+1}=T_{j+1}f_j$. It is not difficult to see that the sequence
$f_j$ converges in $L^1(0,1)$ to a function $f_{\infty}$. Moreover,
the essential oscillation of $f_{\infty}$ over any open set is $2$.
\enddemo

4.5. {\it Linear systems}. 
The examples above can be used to answer  open questions 
(raised in \cite{GS 85})
concerning solutions of linear $2\times 2$ systems of the form
\begin{equation}\label{lineq}
\partial_{\alpha}\aabij(x)\partial_{\beta}v_j=0\, , \qquad i=1,2
\end{equation}
where the coefficients are in $L^{\infty}$ and satisfy  the strong
 Legendre-Hadamard
condition
$$
\aabij(x)\xi_{\alpha}\xi_{\beta}\bar u^i\bar u^j\ge\nu |\xi|^2|\bar u|^2
$$
for each $\xi, \,\bar u\in\rd$ and almost every $x$. (As usual, $\nu>0$.)
In what follows we will write the system~(\ref{lineq}) as 
$\diw A(x) \nabla v=0$.

There is a well known procedure for passing from solutions of nonlinear equations 
to solutions of linear equations with measurable coefficients (see e.g.\ \cite{Mo 66}).
We will use it to construct our examples. These examples will be based on the
following proposition.

\spn{4.1}\proclaim{Proposition}\label{compactsup}
There exist  a smooth strictly quasiconvex function 
$F\colon\mdd\to\RR$ with 
uniformly bounded $D^2F$  and a nontrivial Lipschitz function 
$u\colon\rd\to\rd$ which vanishes for $|x|>1$ and satisfies {\rm (}\/weakly\/{\rm )} the equation 
$\diw DF(\nabla u)=0$ is $\rd$.
\endproclaim

{\it Proof.} We will use the notation introduced earlier in Section~\ref{applications}.
We note that the function $F_1$ from Lemma~\ref{fff} satisfies
$DF_1(0)=0$ and therefore the zero matrix belongs to the set
$K_1\subset\mcd$ corresponding to $F_1$. Thus, we see  that the function
$F_0$ in Theorem~\ref{solutions} can be taken so that $DF_0(0)=0$.
Hence the set $K$ corresponding to $F=F_0$ in Theorem~\ref{solutions}
can be taken so that it contains the zero matrix. We know that there are
nontrivial solutions of $\nabla w\in K$ a.e.\  in $\Om$ which vanish
at $\partial\Omega$. Extending $w$ by zero outside $\Omega$, we get solutions
with the required properties.

\spn{4.2}\proclaim{Proposition}\label{linearprop}
There exist $L^\infty$\/{\rm -}\/coefficients $A(x)$ defined in $\rd$ which satisfy the
strong Legendre\/{\rm -}\/Hadamard condition such that weak solutions of the linear system
$\diw A(x) \nabla v = 0$ exhibit the following behavior.
\begin{itemize}
\ritem{(i)} There exists a compactly supported solution $v$ belonging to
the Sobolev space $W^{1,2}$ but not to $W^{1,2+\delta}$ for any $\delta>0$.
\ritem{(ii)} There exists a sequence $v_j,\,j=1,2,\dots$ of Lipschitz solutions which are
supported in $\{x, |x|<1\}${\rm ,} and converge to zero weakly but not strongly in $W^{1,2}$.
\end{itemize}

\endproclaim

{\it Proof.}
Let $F$ and $u$ be as in Proposition~\ref{compactsup} and let
$$
\tilde A(x)=\int_0^1 D^2F(t\nabla u(x))\,dt.
$$
Since $F$ is smooth and $|D^2F|\le c$, $\tilde A(x)$ is a well-defined
$L^\infty$-function. Since $F$ is strongly quasiconvex, it is also strongly
rank-one convex, and therefore $\tilde A(X)$ satisfies the Legendre-Hadamard
condition. Moreover, we have 
$$
\diw \tilde A(x) \nabla u = 
\diw (DF(\nabla u(x))-DF(0))= 0 \quad\mbox{in $\rd$}
$$
in the weak sense.

Let us consider a sequence $B_{a_j,r_j}\subset\{x\in\rd,\,|x|<1\}$ of 
mutually disjoint balls centered at
$a_j$ with radius $r_j>0$ so that $a_j\to 0$ in $\rd$ and $r_j\to 0$.
We let
\begin{eqnarray*}
A(x) & = & D^2F(0) + \sum_{j=1}^\infty\left(\tilde A(r_j^{-1}(x-a_j))-D^2F(0)\right) 
\qquad \mbox{and}\\
v_j(x) & = & u(r_j^{-1}(x-a_j))\,, \qquad j=1,2,\dots\; .
\end{eqnarray*}
The coefficients $A(x)$ are again bounded and satisfy the strong Legendre-Hadamard
condition. We also have $\,\,\diw A(x) \nabla v_j=0, \quad j=1,2,\dots\  $.
The sequence $v_1,v_2,\dots$ gives (ii). To obtain (i), we consider a sequence
$c_1,c_2,\dots$ satisfying $\sum_{j=1}^\infty c_j^2<\infty$ and 
$\sum_{j=1}^\infty c_j^{2+\delta}=\infty$ for each $\delta>0$.
Then $v=\sum_{j=1}^\infty c_j v_j$ has the required properties.

\AuthorRefNames [9999999]


\begin{references}

\bibitem{AF 87} \name{E.~Acerbi} and \name{N.~Fusco}, A regularity theorem for 
minimizers of quasiconvex integrals, {\it Arch.\ Rational Mech.\ Anal\/}.\  {\bf 99}
(1987), 261--281. 

\bibitem{AH 86}   
\name{R.~Aumann} and \name{S.~Hart},
Bi-convexity and bi-martingales,                          
{\it Israel J.\ Math\/}.\ {\bf 54} (1986), 159--180.

\bibitem{Ba 80} \name{J.\ M.\ Ball}, Strict convexity, strong ellipticity and
regularity 
in the calculus of variations, {\it Math.\ Proc.\ Cambridge Philos.\  Soc\/}.\ 
{\bf 87}
(1980), 501--513. 

\bibitem{Ba 90}
\bibline,  Sets of gradients with no rank-one connections,
{\it J.\ Math.\ Pures Appl\/}.\  (9) {\bf 69} (1990),  241--259. 

\bibitem{BJ 87}  
\name{J.\ M.~Ball} and \name{R.\ D.~James}, 
Fine phase mixtures as minimizers of energy, 
{\it Arch.\ Rational Mech.\ Anal\/}.\ {\bf 100} (1987), 13--52.


\bibitem{BFJK 94} \name{K.~Bhattacharya, N.~Firoozye, R.\ D.~James}, and 
\name{R.\ V.~Kohn}, Restrictions on microstructure, 
{\it  Proc.\ Roy.\ Soc.\ Edinburgh A} {\bf  124} (1994), 843--878.


\bibitem{CT 93}
\name{E.~Casadio-Tarabusi}, 
An algebraic characterization of quasi-convex functions,
{\it Ricerche Mat\/}.\ {\bf 42} (1993), 11--24.

\bibitem{CK 88}
\name{M.~Chipot} and \name{D.~Kinderlehrer},
Equilibrium configurations of crystals,
{\it Arch.\ Rational Mech.\ Anal\/}.\ {\bf 103} (1988), 237--277.



\bibitem{DM 97} \name{B.~Dacorogna} and \name{P.~Marcellini}, General existence theorems for 
Hamilton-Jacobi equations in the scalar and vectorial cases, 
{\it Acta Math\/}.\ 
{\bf 178} (1997), 1--37.

\bibitem{DG 68} \name{E.~De Giorgi}, Un esempio di estremali discontinue per
un problema variazionale di tipo ellittico, {\it Boll.\ Un. Mat.\ Ital\/}.\  {\bf 4} (1968),
135--137.

\bibitem{Ev 86} \name{L.\ C.~Evans}, Quasiconvexity and partial regularity in the
calculus 
of variations, {\it Arch.\ Rational Mech.\ Anal\/}.\ {\bf 95} (1986), 227--252.

\bibitem{Fe 69} \name{H.~Federer}, 
{\it Geometric Measure Theory}, Springer-Verlag, New York, 1969.

\bibitem{GS 85}  \name{M.~Giaquinta} and \name{J.~Sou\v cek},  Caccioppoli's inequality and 
Legendre-Hadamard condition, {\it Math.\ Ann\/}.\  {\bf 270} (1985),  105--107.

\bibitem{GM 68} E.~Giusti and M.~Miranda, Un esempio di
soluzioni discontinue per un problema di minimo relativo ad
un integrale regolare del calcolo delle variazioni, {\it Boll.\ Un.\ Mat.\ Ital\/}.\ 
{\bf 2} (1968), 219--226.


\bibitem{Gr 86} M.\ Gromov, {\it Partial Differential Relations},
Springer-Verlag, New York, 1986. 



\bibitem{HLN 96}
\name{W.~Hao, S.~Leonardi}, and \name{J.~Ne\v{c}as},
An example of irregular solution to a nonlinear
Euler-Lagrange elliptic system with real analytic coefficients. 
{\it Ann.\ Scuola Norm.\ Sup.\  Pisa Cl.\ Sci\/}.\  {\bf 23} (1996), 57--67.


\bibitem{Ku 55} 
\name{N.\ H.~Kuiper}, On $C^1$-isometric embeddings.\ I, 
{\it Nederl.\ Akad.\ Wetensch.\ Proc\/}.\ {\bf A 58} (1955), 545--556.



\bibitem{Ma 85}
\name{P.~Marcellini}, Approximation of quasiconvex functions, and lower semicontinuity
of multiple integrals, {\it Manuscripta Math\/}.\ {\bf 51} (1985), 1--28.



\bibitem{MP 98}
\name{J.~Matou\v{s}ek} and \name{P.~Plech\'{a}\v{c}},
On functional  separately convex hulls, 
{\it  Discrete Comput.\ Geom\/}.\  {\bf 19} (1998), 105--130.










\bibitem{Mo 66} \name{C.\ B.~Morrey}, 
{\it Multiple Integrals in the Calculus of
Variations}, Springer-Verlag, New York, 1966.


\bibitem{MS 96}
\name{S.~M\"uller} and \name{V.~\v{S}ver\'{a}k},
Attainment results for the two-well problem by convex integration, in
{\it Geometric Analysis and the Calculus of Variations} (J.\ Jost, ed.),
Internat.\  Press, Cambridge, MA, 1996, 239--251.



\bibitem{Na 54} \name{J.~Nash}, $C^1$ isometric imbeddings, 
{\it Ann.\ of Math\/}.\ {\bf 60} (1954), 383--396.



\bibitem{NM 91} \name{V.~Nesi} and \name{G.\ W.~Milton}, 
 Polycrystalline configurations that maximize 
electrical resistivity, {\it J.\ Mech.\ Phys.\ Solids} {\bf 39} (1991),
525--542.





\bibitem{Pe 93} \name{P.~Pedregal}, Laminates and microstructure, 
{\it European J.\ Appl.\ Math\/}.\ {\bf 4} (1993), 121--149.

\bibitem{Sch 74} \name{V.~Scheffer}, Regularity and irregularity
of solutions to nonlinear second order elliptic systems of partial
differential equations and inequalities,  Dissertation,
Princeton University, 1974 (unpublished).

\bibitem{Sv 90}
\name{V.~\v{S}ver\'{a}k},
Examples of rank-one convex functions,
{\it Proc.\ Roy.\ Soc.\ Edinburgh} {\bf A 114} (1990), 237--242.



\bibitem{Sv 92a} 
\bibline, 
Rank-one convexity does not imply quasiconvexity,
{\it Proc.\ Roy.\ Soc.\break Edinburgh} {\bf A 120} (1992), 185--189.

\bibitem{Sv 92b} 
\bibline, 
New examples of quasiconvex functions,
{\it Arch.\ Rational Mech.\ Anal\/}.\  {\bf 119} (1992), 293--300. 






\bibitem{Ta 93} \name{L.~Tartar}, Some remarks on separately convex
functions,
in: {\it Microstructure and Phase Transitions},
{\it IMA Vol.\ Math.\ Appl\/}.\  {\bf 54} (D.~Kinderlehrer, R.\ D.~James,
M.~Luskin and J.\ L.~Ericksen, eds.), Springer-Verlag, New York (1993),  191--204. 
 



\end{references}
 \end{document}